\documentclass{amsart}

\usepackage[utf8]{inputenc}
\usepackage{amsthm, amssymb}
\usepackage[colorlinks]{hyperref}
\usepackage[colorinlistoftodos]{todonotes}
\usepackage{pgf}
\usepackage{tikz}
\usepackage{tikz-cd}
\usepackage[numbers]{natbib}
\usetikzlibrary{positioning,shapes,shadows,arrows}

\usepackage{ytableau}

\newtheorem{prop}{Proposition}
\newtheorem{thm}[prop]{Theorem}

\newtheorem{lemma}[prop]{Lemma}
\newtheorem{cor}[prop]{Corollary}

\theoremstyle{remark}
\newtheorem{exa}[prop]{Example}
\newtheorem{rem}[prop]{Remark}
\newtheorem{defn}[prop]{Definition}
\numberwithin{prop}{section}
\numberwithin{equation}{section}

\usepackage[margin=0.75in]{geometry}
\usepackage{amssymb}

\newcommand{\fG}{\mathfrak{G}}
\newcommand{\fGb}{\fG^{(\beta)}}
\newcommand{\fL}{\mathfrak{L}}
\newcommand{\fLb}{\fL^{(\beta)}}

\newcommand{\Gb}{G^{(\beta)}}

\newcommand{\Hequiv}{\equiv_H}

\newcommand{\Inc}{\mathrm{Inc}}
\newcommand{\sInc}{\mathrm{sInc}}

\newcommand{\key}{\mathrm{key}}
\newcommand{\capn}{\mathrm{cap}_n}

\newcommand{\rev}{\mathrm{rev}}
\newcommand{\revKjdt}{\mathrm{revKjdt}}

\newcommand{\RSVT}{\mathrm{RSVT}}

\newcommand{\shape}{\mathrm{shape}}

\newcommand{\word}{\mathrm{word}}
\newcommand{\wt}{\mathrm{wt}}
\newcommand{\Z}{\mathbb{Z}}
\newcommand{\C}{\mathcal{C}}

\newcommand{\T}{\mathcal{T}}
\newcommand{\tl}{\triangleleft}

\definecolor{darkblue}{rgb}{0.0,0,0.7} 
\definecolor{darkred}{rgb}{0.7,0,0} 
\definecolor{darkgreen}{rgb}{0, .6, 0} 
\newcommand{\definition}[1]{{\color{darkred}\emph{#1}}} 

\title{Lascoux expansion of the product of a Lascoux and a stable Grothendieck}
\author[G.~Orelowitz]{Gidon Orelowitz}
\address[G.~Orelowitz]{}
\email{gidon.orelowitz@gmail.com}

\author[T.~Yu]{Tianyi Yu}
\address[T. Yu]{Department of Mathematics, UC San Diego, La Jolla, CA 92093, U.S.A.}
\email{tiy059@ucsd.edu}

\begin{document}
\maketitle

\begin{abstract}

This paper gives a tableau formula for expanding the product of a
Lascoux polynomial and a stable Grothendieck polynomial into Lascoux polynomials.
Lascoux and stable Grothendieck polynomials are inhomogeneous analogues 
of key polynomials and Stanley symmetric functions, respectively. 
Our formula refines the
K-theoretic Littlewood-Richardson rule of Buch 
and extends the key expansion
of key times Schur established 
by Haglund, Luoto, Mason, and van Willigenburg.
Our proof is combinatorial, relying heavily on a novel row insertion algorithm 
of Huang, Shimozono and Yu.

\end{abstract}

\section{Introduction}
\label{S: Intro}
Fix $n\in \mathbb{Z}_{>0}$ throughout this paper.
We consider two families of 
inhomogeneous polynomials.
Both families have positive
integer coefficients 
and involve variables $\beta, x_1, \cdots, x_n$.
\begin{itemize}
\item 
Introduced by Lascoux~\cite{Las}, 
the \definition{Lascoux polynomial} $\fLb_\alpha$
is indexed by a \definition{
weak composition} $\alpha$,
a sequence of $n$ non-negative integers.
\item 
Introduced by Fomin and Kirillov~\cite{FK},
the \definition{stable Grothendieck polynomial} $G^{(\beta)}_w$
is indexed by a permutation $w \in S_+$,
the set of permutations of $\mathbb{Z}_{>0}$
where only finitely many numbers are permuted.
The polynomial $G^{(\beta)}_w$
involves $\beta$ and $x_1, x_2, \cdots$.
We let $G^{(\beta)}_w(x_1, \cdots, x_n)$
be the polynomial obtained by setting $x_{n+1} = \cdots = 0$ in $G^{(\beta)}_w$.
\end{itemize}

We give a tableau formula
that expands $\fLb_\alpha \times G^{(\beta)}_w(x_1, \cdots, x_n)$ into Lascoux polynomials,
where the coefficients are positive integers
multiplied by a power of $\beta$.
Our expansion simultaneously extends the following two 
results.
Both results are generalizations of the famous Littlewood-Richardson rule that gives the Schur expansion of the product of two Schur polynomials $s_\lambda$.

\begin{itemize}
\item[(1)] 
Introduced by Demazure~\cite{Dem},
the key polynomials $\kappa_\alpha$,
are characters of
the Borel subgroup $B$ of upper triangular matrices in $GL_n$.
They can be viewed as non-symmetric generalization 
of Schur polynomials $s_\lambda$:
when $\alpha = (\alpha_1, \cdots, \alpha_n)$ is weakly increasing, 
$\kappa_\alpha$ agrees with $s_{(\alpha_n, \cdots, \alpha_1)}(x_1, \cdots, x_n)$~\cite{Dem}.
Haglund, Luoto, Mason, and van Willigenburg~\cite{HLSV}
established a non-symmetric refinement of 
the Littlewood-Richardson rule:
They expanded
$\kappa_\alpha \times s_\lambda(x_1, \cdots, x_n)$
into key polynomials
using skyline fillings.
\item[(2)] 
In enumerative geometry, 
Schur polynomials represent Schubert classes 
in the cohomology ring of the Grassmannian.  
When $w$ is a Grassmannian permutation,
the stable 
Grothendieck polynomial $G^{(\beta)}_w$ 
is the connective K-theoretic analogue of Schur polynomials~\cite{LS:Groth, Hud}.
Combinatorially, this means if we set 
$\beta$ to be $0$ in $G^{(\beta)}_w$
when $w$ is Grassmannian, 
we get a Schur polynomial. 
Buch~\cite{B} established the connective
K-theoretic Littlewood-Richardson rule: 
a tableau formula that computes the
coefficient of $G^{(\beta)}_w$ in $G^{(\beta)}_u \times G^{(\beta)}_v$
for Grassmannian permutations $u, v, w$.
\end{itemize}

Lascoux polynomials generalize both 
key polynomials and $G^{(\beta)}_w$
for Grassmannian $w$.
If we set 
$\beta$ to be $0$ in $\fLb_\alpha$,
we get the key polynomial $\kappa_\alpha$~\cite{Las}.
If $\alpha$ is weakly increasing, 
$\fLb_\alpha$ agrees with $G^{(\beta)}_w(x_1, \cdots, x_n)$
for some Grassmannian $w$~\cite{BSW}.
Therefore, the expansion we are 
studying is an inhomogeneous
analogue of (1) 
and a non-symmetric extension of (2). 
 
When $w$ is Grassmannian, 
Monical~\cite{Mon} 
described a conjectural rule for this expansion
involving genomic semistandard skyline fillings.
Monical's conjecture is still open.
In this paper, we establish a different
combinatorial rule for arbitrary $w \in S_+$.
Our rule involves
\definition{increasing tableaux}, fillings of Young diagrams with positive integers 
such that each row and column is strictly increasing. 
Our rule applies a sequence of 
three operators to increasing tableaux:
$K_-(\cdot), \capn(\cdot)$ and $\wt(\cdot)$.

For a tableau $T$,
we use $T_i$ to denote the set of numbers
that appear in column $i$ of $T$.
We say a tableau $T$ is a \definition{key}
if each of its column strictly increasing and $T_1 \supseteq T_2 \supseteq \cdots$.
Each increasing tableau $P$
is associated with a key called 
its \definition{left key}, 
denoted as $K_-(P)$. 
In Section~\ref{SS: Kjdt}, we give the usual definition of $K_-(P)$ using
the K-theoretic jeu-de-taquin of Thomas and Yong~\cite{TY}.
In Section~\ref{SS: left key and tl},
we derive a method to compute $K_-(P)$
using the $\tl$ operator.

Take finite sets $S,T \subseteq \mathbb{Z}_{>0}$.
Define $T \tl S$ via the following algorithm.
Go through elements in $S$ from the largest to the smallest. 
For $s$ in $S$,
it picks the largest number in $T$ that is less than $s$
and has not been picked.
If such a number exists, we put it in $T \lhd S$.  
For instance, $\{1,3,4,6,7,9\} \triangleleft \{ 2,3,7,8\}
= \{1,6,7\}$.
In Corollary~\ref{C: leftkeycalc}, 
we show that $K_-(P)_j$ consists of $P_1 \tl P_2 \tl \cdots \tl P_j$
where the evaluation of $\tl$ is from right to left.

 The operator $\capn(\cdot)$ is defined 
on keys with at most $n$ rows.
It replaces numbers larger than $n$ in each column 
by largest integers in $[n]$ that are missing
in that column.
Then it sorts each column to make it increasing. 
For instance, a column consisting of 
$1, 4, 6,8, 9$ will consist of $1,3,4,5,6$
after $\textrm{cap}_6$.

Finally, let $\wt(\cdot)$
be the operator that sends a tableau with entries in 
$[n]$
to a weak composition.
The $i^\textsuperscript{th}$ entry 
is the number of $i$'s in the tableau.

\begin{exa}
\label{E: Main}
Suppose $n = 3$.
Consider the following three increasing tableaux:

$$
\begin{ytableau}
1 & 4 & 6 & 7\\
3\\ 
7
\end{ytableau}
\quad\quad\quad\quad\quad
\begin{ytableau}
1 & 4 & 6 & 7\\
3 & 7\\ 
7
\end{ytableau}
\quad\quad\quad\quad\quad
\begin{ytableau}
1 & 4 & 6 & 7\\
3 & 6 \\ 
6 & 7
\end{ytableau}
$$
After applying $K_-(\cdot)$,
they become the following three keys:

$$
\begin{ytableau}
1 & 3 & 3 & 3\\
3\\ 
7
\end{ytableau}
\quad\quad\quad\quad\quad
\begin{ytableau}
1 & 1 & 3 & 3\\
3 & 3\\ 
7
\end{ytableau}
\quad\quad\quad\quad\quad
\begin{ytableau}
1 & 1 & 3 & 3\\
3 & 3 \\ 
6 & 6
\end{ytableau}
$$

After applying $\textrm{cap}_3(\cdot)$,
they become 

$$
\begin{ytableau}
1 & 3 & 3 & 3\\
2\\ 
3
\end{ytableau}
\quad\quad\quad\quad\quad
\begin{ytableau}
1 & 1 & 3 & 3\\
2 & 3\\ 
3
\end{ytableau}
\quad\quad\quad\quad\quad
\begin{ytableau}
1 & 1 & 3 & 3\\
2 & 2 \\ 
3 & 3
\end{ytableau}
$$

Finally, we apply $\wt(\cdot)$
and get $(1,1,4), (2,1,4)$ and $(2,2,4)$.
\end{exa}

 Each increasing tableau $P$
is associated with a word denoted as $\word(P)$
known as the \definition{reading word}.
It is obtained by reading the entries of $P$ from left
to right, and bottom to top in each column.
The three increasing tableaux in Example~\ref{E: Main}
have reading words $731467$, $7317467$, and $63176467$
respectively. 
Each word $a$ is associated with a permutation
$[a]_H \in S_+$ defined as follows. 
The \definition{0-Hecke monoid} is the quotient of the free monoid of words 
on the alphabet $[n-1]$ by the relations
\begin{align*}
  ii &\Hequiv i \\
  i(i+1)i &\Hequiv (i+1)i(i+1) \\
  ij&\Hequiv ji \qquad\text{for $|i-j|\ge 2$.}
\end{align*}

Let $s_i \in S_+$ be the 
transposition that swaps $i$ and $i+1$.
A \definition{reduced word} of $w \in S_+$
is a word $a_1 a_2 \cdots a_l$ with minimal length 
such that $w = s_{a_1} \cdots a_{a_l}$.
We denote $\ell(w) = l$.
We say $a$ is a \definition{Hecke word} 
of $w$,
denoted as $[a]_H = w$,
if $a \Hequiv b$ and $b$ is a reduced word of $w$.

\begin{exa}
We have $421433 \Hequiv 2143$,
which is a reduced word of $w = s_2 s_1 s_4 s_3$.
Thus $[421433]_H = [2143]_H = w$.
\end{exa}

For a weak composition $\alpha = (\alpha_1, \cdots, \alpha_n)$,
we define $|\alpha| := \alpha_1 + \cdots + \alpha_n$.
For a tableau $P$, we define $|P|$ as 
the number of cells in $P$.
Now we can describe our main result. 

\begin{thm}
\label{thm:G times L}
Let $\alpha$ be a weak composition.
Take an increasing tableau $P_1$
with $\wt(K_-(P_1)) = \alpha$.
Let $N$ be a number such that $N > n$
and $N > \max(P_1)$.
Suppose $w \in S_+$
fixes $1, 2, \cdots, N$.
Then
\begin{align}
\label{EQ: G times L}
\fLb_\alpha \times G^{(\beta)}_w(x_1, \cdots, x_n)
= \sum_{P} \beta^{|P| - \ell(w) - |\alpha|} \fLb_{\wt(\capn(K_-(P)))},
\end{align}
where the sum is over all increasing tableau $P$ which has at most $n$ rows 
and satisfies the following.
\begin{itemize}
\item If we ignore entries of $P$ 
that are larger than $N$,
we get $P_1$.
\item If we ignore all entries smaller than $N$ in $\word(P)$, we get a Hecke word for $w$.
\end{itemize}
\end{thm}

We explain why Theorem \ref{thm:G times L} gives the Lascoux expansion of $\fLb_\alpha \times G_w(x_1, \cdots, x_n)$ 
for arbitrary weak composition $\alpha$ and $w \in S_+$.
First, we can always find increasing tableau $P_1$
with ${\sf wt}(K_-(P_1)) = \alpha$ for any $\alpha$:
If $\alpha = (\alpha_1, \cdots, \alpha_n)$,
we may let $P_1$ be the increasing tableau
whose column $c$ consists of $\{\alpha_i + c -1: \alpha_i \geq c \}$.
Next, a fact is $G^{(\beta)}_w(x_1, \cdots, x_n) = G^{(\beta)}_{1^N \times w}(x_1, \cdots, x_n)$,
where $1^N \times w$ sends $i$ to $w(i-N)+N$ for all $i > N$ and fixes $1, 2, \cdots, N$.
This fact can be deduced from~(\ref{EQ: stable G}).
Thus, we may replace $w$ by $1^N \times w$,
so $w$ satisfies the conditions
in Theorem~\ref{thm:G times L}.

\begin{exa}
Say we would like to expand 
$\fLb_{\alpha} \times 
G^{(\beta)}_{w}(x_1, \cdots, x_n)$ 
into Lascoux polynomials
where $\alpha = (1,0,2)$, 
$n = 3$ 
and $w$ has one-line notation
$321$.
We may let $P_1$ be

$$
P_1 = 
\raisebox{0.1cm}{
\begin{ytableau}
1 & 4\\
3
\end{ytableau}}
$$
so that  $K_-(P_1) = (1,0,2)$.
Then we may pick $N = 5$, 
so that $N > \max(P_1)$
and $N > n$.
We may redefine $w$ 
as the permutation with 
one-line notation 
$12345876$.
This replacement does not
change $G^{(\beta)}_w$.
Then the three increasing 
tableaux in Example~\ref{E: Main}
satisfy the conditions
of Theorem~\ref{thm:G times L}.
They contribute 
$\fLb_{(1,1,4)}$, $\beta\fLb_{(2,1,4)}$
and $\beta^2\fLb_{(2,2,4)}$
to the expansion. 
In Appendix \ref{appendix}, 
we enumerate all increasing tableaux 
that satisfy the conditions for this example. 

\end{exa}

\begin{rem}

After setting $\beta = 0$,
$\fLb_\alpha$
becomes $\kappa_\alpha$
and $\Gb_w(x_1, \cdots, x_n)$
becomes $F_w(x_1, \cdots, x_n)$,
the Stanley symmetric function~\cite{St}
in variables $x_1, \cdots, x_n$.
By setting $\beta = 0$ in~(\ref{EQ: G times L}), we have:
\begin{align}
\kappa_\alpha \times F_w(x_1, \cdots, x_n)
= \sum_{P} \kappa_{\wt(\capn(K_-(P)))},
\end{align}
where the sum is over all 
increasing tableau $P$ 
such that 
\begin{itemize}
\item If we ignore entries of $P$ 
that are larger than $N$,
we get $P_1$.
\item If we ignore all entries smaller than $N$ in $\word(P)$, we get a \textit{reduced word} for $w$.
\end{itemize}

When $w$ is Grassmannian, 
$F_w(x_1, \cdots, x_n)$ is a Schur polynomial in $x_1, \cdots, x_n$.
Thus, our result restricts to a 
tableau formula for the key expansion
of $\kappa_\alpha \times s_\lambda(x_1, \cdots, x_n)$.
We do not have a bijection
between our tableaux
and skyline fillings in~\cite{HLSV}.
\end{rem}

Our approach relies on combinatorial 
formulas of $G^{(\beta)}_w$ and $\fLb_\alpha$.
Fomin and Kirillov established a 
combinatorial formula of $G^{(\beta)}_w$
via compatible pairs, certain pairs of words.
Buciumas, Scrimshaw, and Weber~\cite{BSW}
established a formula for Lascoux polynomials involving reverse set-valued tableaux (RSVT).
Our main tool is an insertion algorithm
developed by Huang, Shimozono and Yu~\cite{HSY}.
This algorithm is a row insertion analogue 
of Hecke column insertion~\cite{BKSTY}.
It gives a bijection between compatible pairs
and the set of $(P, Q)$ 
where $P$ is an increasing tableau
and $Q$ is a RSVT of the same shape.
In Theorem~\ref{thm:insertionbijection}, 
we show this row insertion satisfies
an analogous property of Hecke column insertion
established in~\cite[Theorem 4.2]{SY}.
This property allows us to 
turn the tableau formula of $\fLb_\alpha$
into a compatible pair formula. 
Then both sides of Theorem~\ref{thm:G times L}
can be viewed as a sum of certain compatible pairs.
Finally, we establish bijections between the compatible pairs representing the two sides.

\begin{rem} 

In~\cite{SY}, the authors established Theorem 4.2
for Hecke column insertion 
under the convention that the insertion tableau 
is a decreasing tableau. 
The analogous property would not hold 
for Hecke column insertion of increasing tableaux.
However, we do not have a way to state
Theorem~\ref{thm:G times L} using decreasing tableaux.
Moreover, our main argument would not work for decreasing tableaux (see Remark~\ref{R: Why not Hecke}).
Thus, we chose to use the row insertion of Huang, Shimozono and Yu~\cite{HSY} on increasing tableaux.
\end{rem}

Lastly, 
a by-product of Theorem~\ref{thm:insertionbijection}
is the Lascoux expansion of Grothendieck polynomials. 
This expansion was first conjectured by Reiner and Yong~\cite{RY}
and proved by Shimozono and Yu~\cite{SY}.
The proof in~\cite{SY} uses 
Hecke column insertion on decreasing tableaux
while our proof uses 
the row insertion on increasing tableaux.

The rest of the paper is organized as follows. 
In \S\ref{S: Background}, 
we cover some necessary background. 
In \S\ref{S: Proof},
we state Theorem~\ref{thm:insertionbijection}
and use it to prove Theorem~\ref{thm:G times L}.
The rest of the paper aims to prove 
Theorem~\ref{thm:insertionbijection}.
In \S\ref{S: triangle},
we introduce and investigate the operator $\lhd$.
In \S\ref{S: left key},
we define the left key of increasing tableaux
and give a simple way to compute it using $\lhd$.
In \S\ref{S: insertion},
we describe and study the insertion algorithm
of Huang, Shimozono and Yu.
Then we prove Theorem~\ref{thm:insertionbijection}.

\section{Background}
\label{S: Background}

In this section, we first introduce 
the two main players of this paper: 
stable Grothendieck polynomials and Lascoux polynomials. 
Instead of providing their usual definitions
in the literature,
we describe combinatorial formulas
to compute them.
Finally, we briefly describe the 
reverse row insertion
and leave the details to \S\ref{S: insertion}.

\subsection{Stable Grothendieck polynomials
and compatible pairs}

The stable Grothendieck polynomial can be computed using 
compatible pairs. 

\begin{defn} \label{D:compatible pair} \cite{BJS}
Let $a = a_1 \cdots a_m$
and $i = i_1 \cdots i_m$
be two words of the same length.
We say $(a,i)$
is a \definition{compatible pair} if
\begin{itemize}
\item $i$ is weakly increasing, and 
\item $i_j = i_{j+1}$ implies $a_j > a_{j+1}$.
\end{itemize}
We say the compatible pair $(a,i)$ is \definition{bounded} 
if $i_j \leq a_j$ for all $j \in [m]$.
\end{defn}

Let $\C$ be the set of all compatible pairs.
Let $\C^b$ be the set of bounded compatible pairs.
Let $\C_w$ (resp. $\C_w^b$) be the set of $(a,i) \in C$ (resp. $C^b$)
such that $a$ is a Hecke word of $w$.

\begin{exa}
The pair $(421433, 111224)$ is compatible and thus in $\C$.
This pair is not bounded since the last number in the first word
is smaller than the last number in the second word. 
On the other hand, $(421433, 111223)$ is bounded. 
\end{exa}

For a word $i$, let $\ell(i)$ be its length
and let $\wt(i)$ be a sequence where the $j\textsuperscript{th}$ entry is the 
number of times $j$ appears in $i$.
For a sequence of numbers $(c_1, c_2, \cdots)$
with only finitely many non-zero entries,
we use $x^{(c_1, c_2, \cdots)}$ to denote the 
monomial where the power of $x_i$ is $c_i$,
\begin{defn}~\cite{FK}
\label{D: Gro} 
The \definition{Grothendieck polynomial} $\fGb_w$ and 
the \definition{stable Grothendieck function} $G^{(\beta)}_w$
can be defined as:
\begin{align}
\label{EQ: GrotDef}
\fGb_w = \sum_{(a,i)\in \C^b_{w^{-1}}} \beta^{\ell(a) - \ell(w)} x^{\wt(i)}\\
\Gb_w = \sum_{(a,i)\in \C_{w^{-1}}} \beta^{\ell(a) - \ell(w)} x^{\wt(i)}.
\end{align}
\end{defn}

\begin{exa}
Consider $w$ with one-line notation $31524$.
Then $w^{-1}$ has one-line notation $24153$.
We know $[421433]_H = w^{-1}$ and $(421433, 111223)$ is bounded,
so this pair is in $C^b_{w^{-1}}$.
It would contribute $\beta^2 x_1^3 x_2^2 x_3$
to $\fGb_w$ and $\Gb_w$.
\end{exa}

Notice that the stable Grothendieck function generally 
involves an infinite set of variables $x_1, x_2, \cdots$ and has infinitely many terms. 
One main player of this paper is $G^{(\beta)}_w(x_1, \cdots, x_n)$, which is obtained from $G^{(\beta)}_w$
by setting $x_{n + 1} = \cdots = 0$.
Clearly, 
\begin{align}
\label{EQ: stable G}
\Gb_w(x_1, \cdots, x_n) = \sum_{(a,i)\in \C_{w^{-1}}^{\leq n}} \beta^{\ell(a) - \ell(w)} x^{\wt(i)},   
\end{align}
where
$\C_{w}^{\leq n} := \{(a, i) \in \C_w: \textrm{entries of $i$ are at most $n$}\}.$

\subsection{Lascoux polynomials and tableaux}
A \definition{weak composition} is a sequence of $n$ non-negative integers.
For a weak composition $\alpha$,
we use $\alpha_i$ to denote
its $i\textsuperscript{th}$
entry. 
We also define $x^\alpha$ as the monomial
$x_1^{\alpha_1}\cdots x_n^{\alpha_n}$
and define $|\alpha| := \sum_{i = 1}^n\alpha_i$.

For a weak composition $\alpha$, Lascoux~\cite{Las} defined 
the \definition{Lascoux polynomial} $\fLb_\alpha \in \Z_{\geq 0}[\beta][x_1, \cdots, x_n]$.
We define $\fLb_\alpha$
using a tableaux formula
of Buciumas, Scrimshaw and Weber ~\cite{BSW}.
It generalizes a classical tableau formula
of key polynomials found by Lascoux and Sch{\"u}tzenberger~\cite{LS:key}.

A partition is a weakly decreasing sequence of positive numbers. 
The \definition{Young diagram} of a partition $\lambda = (\lambda_1, \dots, \lambda_m)$ is 
the set $\{(r,c): c \leq \lambda_r\}$.
We represent a Young diagram by drawing a cell in row $r$ column $c$
for each $(r,c)$ in the set
under the English convention: 
Row $1$ is the topmost row
and column 1 is the leftmost column.
A \definition{tableau} is a filling of a Young diagram.
For a tableau $T$, we use $T(r,c)$ to denote 
its filling in the cell $(r,c)$.
The \definition{shape} of $T$, denoted as $\shape(T)$,
is the underlying Young diagram of $T$.
We say $(r,c)$ is a cell in $T$ if $(r,c)$
is in $\shape(T)$.

In this paper, usually, we fill each cell by $\mathbb{Z}_{>0}$
or subsets of $\mathbb{Z}_{>0}$.
When $T$ is a tableau filled by $[n]$,
we let $\wt(T)$ be the weak composition 
whose $i\textsuperscript{th}$
entry is the number of $(r,c)$ in $T$ 
such that $T(r,c) = i$.
When $T$ is a tableau filled by subsets of $[n]$,
we let $\wt(T)$ be the weak composition 
whose $i\textsuperscript{th}$
entry is the number of $(r,c)$ in $T$ 
such that $i \in T(r,c)$.

We define Lascoux polynomials
using tableaux.
A \definition{reverse semi-standard Young tableau (RSSYT)} is a 
filling of the Young diagram with $[n]$.
We require every row (resp. column) to be weakly 
(resp. strictly) decreasing from left to right 
(resp. top to bottom). 
A \definition{reverse set-valued tableau} (RSVT) is a filling of the Young diagram
with non-empty subsets of $[n]$.
Moreover, for two horizontally adjacent cells $(r,c)$ and $(r,c+1)$,
we require $\min(T(r,c)) \geq \max(T(r,c+1))$.
For two vertically adjacent cells $(r,c)$ and $(r+1,c)$,
we require $\min(T(r,c)) > \max(T(r+1,c))$.
For an RSVT $T$,
let $L(T)$ be the tableau obtained by keeping only
the largest number in each entry of $T$.
Clearly, $L(\cdot)$ is a map from RSVTs to RSSYTs. 

\begin{exa}
\label{example:RSVT+RSSYT}
For the following RSVT $T$,
we compute $L(T)$.
\[
T = \raisebox{0.4cm}{
\begin{ytableau}
6 & 53 & 3 & 321\\
4 & 21 & 1\\
32
\end{ytableau}},
\quad\quad\quad\quad
L(T) = \raisebox{0.4cm}{
\begin{ytableau}
6 & 5 & 3 & 3\\
4 & 2 & 1\\
3
\end{ytableau}}.
\]
When writing an entry in a RSVT,
we simply list its elements in decreasing order. 
For instance, the ``$6$'' in cell $(1,1)$
represents $T(1,1) = \{6\}$
and ``$321$'' in cell $(1,4)$
represents $T(1,4) = \{3,2,1\}$.
\end{exa}

Lascoux polynomials
can be written as sums over RSVTs. 
To describe which RSVTs can appear in a sum,
we need a few more definitions.
Recall a \definition{key} is
a tableau where
each column is increasing
and numbers in one column also appear in the 
column to the left.
There is a natural bijection, denoted as $\key(\cdot)$, 
that sends a weak composition to a key
with entries in $[n]$.
For a weak composition $\alpha$,
$\key(\alpha)$ is the unique key 
whose column $c$ consists of $\{i: \alpha_i \geq c\}$.
Its inverse is just $\wt(\cdot)$.
\begin{exa}
Let $\alpha = (2,1,4,0,2)$.Then

\[
\key(\alpha) = \raisebox{0.6cm}{
\begin{ytableau}
1 & 1 & 3 & 3\\
2 & 3\\
3 & 5\\
5 \\
\end{ytableau}}.
\]
\end{exa}

Each RSSYT $T$ is associated 
with a key called its \definition{left key},
denoted as $K_-(T)$.
The left key is originally defined via moves
on RSSYT known as jeu-de-taquin.
We present one simple method to compute
the left key, which is a reformulation of 
Willis' algorithm~\cite{Wil},
using an operation on sets. 
\begin{defn}[{\cite[Definition 3.11]{SY}}]
\label{D: Triangle}
Take finite sets $S, T \subseteq \Z_{>0}$.
Define $T \trianglerighteq S$ be the set
computed as follows. 
Iterate through elements of $S$
from the smallest to the largest. 
For each $s \in S$,
it picks the smallest $t \in T$ such that $t \geq s$
and $t$ has not been picked.
Then $T \trianglerighteq S$ is the set
of all picked $t \in T$.   
\end{defn}
For instance, $\{ 1,4,5,6,7\} \trianglerighteq \{ 2,3,7\} = \{4,5,7\}$, 
since $2$ picks $4$, $3$ picks $5$ and $7$ picks $7$.

\begin{defn}
\label{D: left key of RSSYT}
For a RSSYT $T$, its left key $K_-(T)$
is the key whose column $i$ consists of
$$
T_1 \trianglerighteq (T_2 \trianglerighteq(\cdots (T_{i-1} \trianglerighteq T_i)\cdots )).
$$
\end{defn}

\begin{exa}
Let $T$ be the following RSSYT:
$$
T = \raisebox{0.4cm}{
\begin{ytableau}
6 & 5 & 3 & 3\\
4 & 2 & 1\\
3
\end{ytableau}}.
$$
Then column $1$ of $K_-(T)$
consists of $T_1 = \{3,4,6\}$.
Column $2$ of $K_-(T)$
consists of 
$T_1 \trianglerighteq T_2 = \{3,6\}$.
Column $3$ of $K_-(T)$
consists of 
$T_1 \trianglerighteq (T_2 \trianglerighteq T_3) = \{3,6\}$.
Column $4$ of $K_-(T)$
consists of 
$T_1 \trianglerighteq (T_2 \trianglerighteq (T_3 \trianglerighteq T_4)) = \{6\}$. Thus, 
$$
K_-(T) = \raisebox{0.4cm}{
\begin{ytableau}
3 & 3 & 3 & 6\\
4 & 6 & 6\\
6
\end{ytableau}}.
$$
\end{exa}

\begin{defn}
We define the Lascoux polynomial $\fLb_\alpha$
using a tableau formula~\cite{BSW}:
\begin{equation}
\label{EQ: RSVT rule for Lascoux}
\fLb_\alpha = \sum_{
\substack{\textrm{RSVT } T\\ K_-(L(T)) \leq \key(\alpha)}}  \beta^{|\wt(T)| - |\alpha|}x^{\wt(T)}.
\end{equation}
\end{defn}

\begin{exa}
When $\alpha = (0,2,1)$,
the following are all 
the RSVT that contribute
to $\fLb_\alpha$:
$$
\begin{ytableau}
3  & 2 \cr
2
\end{ytableau}\,,
\quad  
\begin{ytableau}
2  & 2 \cr
1
\end{ytableau}\,,
\quad 
\begin{ytableau}
3  & 1 \cr
2
\end{ytableau}\,,
\quad 
\begin{ytableau}
2  & 1 \cr
1
\end{ytableau}\,,
\quad 
\begin{ytableau}
3  & 1 \cr
1
\end{ytableau}
$$
$$
\begin{ytableau}
3  & 2 \cr
21
\end{ytableau}\,,
\quad  
\begin{ytableau}
3  & 21 \cr
2
\end{ytableau}\,,
\quad 
\begin{ytableau}
2  & 21 \cr
1
\end{ytableau}\,,
\quad 
\begin{ytableau}
32  & 1 \cr
1
\end{ytableau}\,,
\quad 
\begin{ytableau}
3  & 1 \cr
21
\end{ytableau}\,,
\quad 
\begin{ytableau}
3  & 21 \cr
21
\end{ytableau}\,.
$$
Thus, 
$
\fLb_\alpha = 
x_2^2x_3 + x_1x_2^2 + x_1x_2x_3 + x_1^2x_2 + x_1^2x_3 + \beta(2x_1x_2^2x_3 + x_1^2x_2^2 + 2x_1^2x_2x_3)
+ \beta^2 x_1^2x_2^2x_3$.
\end{exa}

\subsection{Row insertion algorithm}
The main tool of this paper is an insertion algorithm of Huang, Shimozono and Yu~\cite{HSY}.
We adopt a slightly different convention: 
They described the insertion algorithm on decreasing tableaux
but we describe it on increasing tableaux.

Let $\T$ be the set of pairs of tableaux $(P, Q)$,
where $P$ is increasing, $Q$ is RSVT and 
$\shape(P) = \shape(Q)$. 
For $w \in S_+$,
let $\T_w$ be the the set $\{(P, Q) \in \T: [\rev(\word(T))]_H = w\}$
where $\rev(\cdot)$ reverses a word.

In~\cite{HSY},
Huang, Shimozono and Yu described a map 
$\Psi: \T \rightarrow \C$
which relies on a row insertion algorithm. 
We postpone the technical
definition of this map to 
\S\ref{S: insertion}.
For now, we just need the following result.
\begin{thm}\cite[Corollary 5.9]{HSY}
The map $\Psi$ is a bijection from
$\T_w$ to $\C_w$
that preserves the weight of the second entry.
In other words, $\Psi$ is a bijection from $\T$ to $\C$
such that if $\Psi(P,Q) = (a,i)$,
then $[\rev(\word(P))]_H \Hequiv [a]_H$
and $\wt(Q) = \wt(i)$.
\end{thm}

We also need the following statement
which we will prove in Section \ref{SS: reverse insertion}.
\begin{lemma}
\label{L: restriction to smaller number}
Take $(A,I) \in \C$
and let $N$ be an integer that
does not appear in $A$.
Let $a$ be the word obtained
by removing all numbers larger than $N$ in $A$.
Suppose $\Psi^{-1}(A,I) = (P,Q)$
and $\Psi^{-1}(a,i) = (p,q)$
for some $i$ such that $(a,i) \in \C$.
Then after removing all numbers
larger than $N$ in $P$,
we obtain $p$.
\end{lemma}

\section{Proof of Theorem~\ref{thm:G times L}}
\label{S: Proof}
Our proof approach requires a technical lemma
which describes the preimage of $\C^b_w$ under $\Psi$.
Let $\T^b$ be the set of $(P,Q) \in \T$
such that $K_-(P) \geq K_-(Q)$
where the comparison is entry-wise. 
Define $\T^b_w$ as the intersection of $\T^b$
and $\T_w$.

\begin{thm}
\label{thm:insertionbijection}
 The map $\Psi$ is a bijection from $\T^b$ to $\C^b$,
 so it is a bijection from $\T^b_w$ to $\C^b_w$
 for any $w \in S_+$.
\end{thm}

In this section, 
we use Theorem~\ref{thm:insertionbijection} to prove
Theorem~\ref{thm:G times L}.
The proof of Theorem~\ref{thm:insertionbijection}
will be the goal of all remaining sections. 

First, we study  
the $\capn(\cdot)$ operator defined in \S\ref{S: Intro}.
Let $T$ be a key
with at most $n$ rows. 
Notice that a number $i$ is in $\capn(T)_c$ if and only if 
$i \in T_c$ or $|T_c \cap (i, \infty)| > n - i$.
Then we have the following property:

\begin{lemma}
\label{L: cap}
Let $T$ be a key with at most $n$ rows. 
Then $\capn(T)$ is also a key. 
Moreover, if $T'$ is another key, 
$\capn(T) \geq T'$ if and only if $T \geq T'$ and $\max(T') \leq n$.
\end{lemma}

\begin{proof}
Assume $i\in \capn(T)_c$ for some $c > 1$.
Then $i \in T_c$ or $|T_c \cap (i, \infty)|> n - i$,
so $i \in T_{c-1}$ or $|T_{c-1} \cap (i, \infty)|> n - i$.
We have $i \in \capn(T)_{c-1}$, so $\capn(T)$ is a key.

Assume $\capn(T) \geq T'$. 
We know entries of $\capn(T)$ are at most $n$,
so $\max(T') \leq n$.
Also, $T \geq \capn(T) \geq T'$.

Now assume $T \geq T'$ and $\max(T') \leq n$.
Consider an arbitrary column $T_c$ with $m$ entries.
If $T_c(i) = \capn(T)_c(i)$,
we clearly have $\capn(T)_c(i) \geq T'_c(i)$.
If $T_c(i) > \capn(T)_c(i)$,
we know $\capn(T)_c(i), \cdots, \capn(T)_c(m)$
are exactly $n -m + i, \cdots, n$.
Since $\max(T') \leq n$,
we have $T_c'(m) \leq n$ so $T_c'(i) \leq n - m + i = \capn(T)_c(i)$.
\end{proof}
 
\begin{cor}
\label{C: cap}
For an key $T$ with at most $n$ rows,
we have  
$$
\sum_{\RSVT \: Q: K_-(L(Q)) \leq T, \max(Q) \leq n} 
\beta^{|\wt(Q)| - |\wt(T)|}x^{\wt(Q)}
= \fLb_{\wt(\capn(T))}. 
$$
\end{cor}
\begin{proof}
Just need to check the $Q$ 
we sum over are precisely all $Q$
such that $K_-(L(Q)) \leq \capn(T)$.
Notice that $\max(Q) \leq n$
is equivalent to $\max(L(Q)) \leq n$,
which is equivalent to $\max(K_-(L(Q))) \leq n$. 
Then the proof is finished by Lemma~\ref{L: cap}.
\end{proof}

Now we prove our main result. 
\begin{proof}[Proof of Theorem~\ref{thm:G times L}]

Observe that the second condition of $P$
is the same as saying $[\word(P)]_H$ 
is a permutation that agrees with $w$
after $N$.

Define
$$
\T_1 = \{(P,Q) \in \T^b: P = P_1\} \textrm{ and } C_1 := \Psi(\T_1).
$$
By Theorem~\ref{thm:insertionbijection},
$C_1 \subset \C^b$.
Let $C_2$ denote $\C_w^{\leq n}$.
Since we assume $w$ fixes $1, \cdots, N$ and $N > n$,
we know $C_2 \subseteq C^b$.
Now by (\ref{EQ: RSVT rule for Lascoux})
and (\ref{EQ: stable G}),
we have
$$
\fLb_\alpha = \sum_{(P, Q) \in \T_1} \beta^{|\wt(Q)| - |\alpha|} x^{\wt(Q)}
= \sum_{(a,i) \in \C_1} \beta^{|\wt(i)| - |\alpha|} x^{\wt(i)}, \textrm{ and}
$$
$$
G^{(\beta)}_w(x_1, \cdots, x_n) 
= \sum_{(b,j) \in \C_2} \beta^{|\wt(j)| - \ell(w)} x^{\wt(j)}.
$$
Let $\T_3$ be the set of $(P, Q) \in \T^b$ such that $P$
satisfies the condition in the Theorem and $\max(Q) \leq n$.
Let $\C_3 = \Psi(\T_3)$,
so $\C_3 \subset \C^b$.
By Corollary~\ref{C: cap},
the right hand side of (\ref{EQ: G times L})
becomes
$$
\sum_{(c,k) \in \C_3} \beta^{|\wt(k)| - \ell(w) - |\alpha|} x^{\wt(k)}.
$$

It remains to build a bijection between $\C_1 \times \C_2$ 
and $\C_3$
such that if $(a, i), (b,j) \mapsto (c,k)$,
then $\wt(i) + \wt(j) = \wt(k)$.
Given $(a,i) \in C_1$
and $(b, j) \in C_2$, numbers in $a$ are smaller than $N$ while numbers in $b$ are
larger than $N$.
There is a unique way to ``shuffle'' them 
and obtain an element of $\C$.
More explicitly, there exists a unique 
$(c,k) = (c_1 \dots c_m, k_1 \cdots k_m)  \in \C$
such that if we let $S = \{ s: c_s <  N\}$ and $T = \{ t: c_t > N\}$
then $$(a,i) = (c_{S(1)} \cdots c_{S(|S|)}, k_{S(1)} \cdots k_{S(|S|)} ) \quad \textrm{ and } \quad (b,j) = (c_{T(1)} \cdots c_{T(|T|)}, k_{T(1)} \cdots k_{T(|T|)} ).$$
Now we check such $(c,k)$ is in $\C_3$.
Since $(a,i), (b,j) \in \C^b$,
so is $(c,k)$.
Let $(P, Q) = \Psi^{-1}(c,j)$.
By Theorem~\ref{thm:insertionbijection},
$\Psi^{-1}(c,j) \in \T^b$.
Clearly, $\max(Q) = \max(k) \leq n$,
which also implies $Q$ and $P$ have at most $n$ rows.
It remains to check 
$P$ satisfies the two conditions in the Theorem:
\begin{itemize}
\item If we only look at numbers in $c$ that are smaller than $N$,
we get $a$.
By Lemma~\ref{L: restriction to smaller number},
if we only look at entries in $P$ that are smaller than $N$,
we get $P_1$.
\item If we only look at numbers in $c$ that are smaller than $N$,
we get $b$, a Hecke word of $w$.
Thus, $[\word(P_1)]_H = [c]_H$ 
is a permutation that agrees with $w$
after $N$.
\end{itemize}

Now we describe the inverse of the bijection above.
Take $(c,k) \in \C_3$.
Let $S = \{ s: c_s <  N\}$ and $T = \{ t: c_t > N\}$
and define $$(a,i) = (c_{S(1)} \cdots c_{S(|S|)}, k_{S(1)} \cdots k_{S(|S|)} ) \quad \textrm{ and } \quad (b,j) = (c_{T(1)} \cdots c_{T(|T|)}, k_{T(1)} \cdots k_{T(|T|)} ).$$
It remains to check $(a,i) \in \C_1$ and 
$(b,j) \in \C_2$.
By $(c,k) \in \C^b$, we know so are $(a,i)$ and $(b,j)$.
By Lemma~\ref{L: restriction to smaller number},
we know $\Psi(a,i) = (P_1, Q)$ for some $Q$,
so $(a,i) \in \C_1$.
Since $b$ is obtained from $c$
by looking at numbers larger than $N$,
$[b]_H$ agree with $[c]_H$ after $N$
and fixes $1, 2, \cdots, N$.
The condition of $\C_3$ guarantees that 
$[c]_H$ agrees with $w$ after $N$.
Thus, $[b]_H = w$, so $(b,j) \in \C_2$. 
\end{proof}

\begin{rem}
\label{R: Why not Hecke}
We explain why our argument works 
on increasing tableaux
but not decreasing tableaux. 
In a decreasing tableau, 
numbers will be decreasing in each column and row. 
If we ignore numbers larger than a given
number, we might not get a tableau of partition shape. 
\end{rem}

\section{The $\lhd$ operator}
\label{S: triangle}

 In this section, we define and study an operator $\lhd$.
In the next section, we will use this operator
to give a simple algorithm that computes $K_-(P)$
of an increasing tableau $P$.

 In the rest of this paper, 
we use $S, T, U$ to denote finite subsets of $\mathbb{Z}_{>0}$.
Given a set $S$, let $S(i)$ be the $i^{th}$ smallest element of $S$ for $i \in [|S|]$.  
Clearly, if $S\subseteq T$, then $S(i) \geq T(i)$ 
for all $i \in [S]$.  
Also, given $x\in S$, we use $S-x$ to denote the set $S\setminus \{x\}$.
Similarly for $x\not\in S$, we use $S \sqcup x$ to denote
the set $S \sqcup \{x\}$.
We evaluate these expressions from left to right.
For instance, 
$S \sqcup x - y = (S \sqcup x) - y$.

Now we define the $\lhd$ operator in two different, but clearly equivalent, ways.

\begin{defn}
Define $T \lhd S$ as follows:
\begin{itemize}
\item Recursive definition: 
If $s\leq t$ for all $s\in S$ and $t\in T$ (including when this is vacuously true, if $S=\emptyset$ or $T = \emptyset$), $T\triangleleft S = \emptyset$.
Otherwise, $\max(S) > \min(T)$.
Let $m = \max(T_{<\max(S)})$ and define $$T\triangleleft S := m \sqcup (T_{<m} \triangleleft (S-\max(S))).$$
\item Non-recursive definition:
We compute $T \triangleleft S$ via the following process.
We initialize $T \lhd S$ as the empty set
and go through elements in $S$ from the largest to the smallest. 
For $s$ in $S$,
it picks the largest number in $T$ that is less than $s$
and has not been picked.
If such a number exists, we put it in $T \lhd S$.  
\end{itemize}  
\end{defn}

We use the recursive definition in proofs by induction and the non-recursive algorithm in other proofs.

\begin{exa}
\label{E: lhd}
We have
$$
\{1,3,4,6,7,9\} \triangleleft \{ 2,3,7,8\}
= \{1,6,7\},\quad\quad
\{1,3,4,6,7,9\} \triangleleft \{ 2,4,7,8\}
= \{1,3,6,7\}.
$$
\end{exa}

Observe as well that the $\triangleleft$ operation is not associative, as 
$$(\{1,2\}\triangleleft \{2,3\}) \triangleleft \{3\} = \{1,2\}\triangleleft \{3\} = \{2\} \quad \textrm{but}\quad \{1,2\}\triangleleft (\{2,3\} \triangleleft \{3\}) = \{1,2\}\triangleleft \{2\} = \{1\}.$$  Because of this, whenever there is ambiguity of the order of evaluation we assume that the expression is evaluated from right-to-left, so $\{1,2\}\triangleleft \{2,3\} \triangleleft \{3\}  = \{1\}$.  

We now study the $\lhd$ operator
from various aspects 
and develop some technical lemmas
that will be used in future sections.

\subsection{Understanding $(T - x) \lhd S$}

For $x \in T$, we would like to understand 
how $(T-x)\triangleleft S$ differs
from $T \lhd S$.
We start with an easy special case:

\begin{lemma}
\label{lemma:singlereduction}
If $x\in T\setminus (T\triangleleft S)$, then 
$(T-x)\triangleleft S = T\triangleleft S$.
\end{lemma}
\begin{proof}
Prove by induction on $|S|$.  If $|S| = 0$,  we are done.  Now assume that the statement is proved whenever $|S| = k-1$, and let $|S|=k$.  If $\max(S) \leq \min(T)$, then $\max(S) \leq \min(T-x)$ as well, and so $(T-x)\triangleleft S = T\triangleleft S = \emptyset$ and we are done.  Otherwise, let $m = \max(T_{<\max(S)})$.  By definition, $m\in T\triangleleft S$, so $m\not=x$ and $m = \max((T-x)_{<\max(S)})$ as well.  
Notice that $(T-x)_{< m}$
is either $T_{< m}$ or $T_{< m} - x$.
In either case, we have 
$(T-x)_{< m}\triangleleft (S-\max(S)) = 
T_{< m}\triangleleft (S-\max(S))$ by the inductive hypothesis.
Thus, 
$$
(T-x) \triangleleft S 
= m \sqcup ((T-x)_{<m} \triangleleft (S-\max(S)))
= m \sqcup (T_{<m} \triangleleft (S-\max(S))) = (T-x)\triangleleft S.
\qedhere
$$
\end{proof}

This Lemma allows us to
remove certain elements from $T$
without changing $T$.

\begin{cor}
\label{lemma:Reduction (subset)}
Let $T' \subseteq T$ such that $T \tl S \subseteq T'$.
Then $T \tl S = T' \tl S$.
\end{cor}
\begin{proof}
This follows from repeated application of Lemma \ref{lemma:singlereduction}.
\end{proof}

Now we analyze the harder case: $x\in T\triangleleft S$.

\begin{lemma}
\label{lemma:singlereduction3}
Take $x\in T\triangleleft S$.  If $T_{<x}\setminus (T\triangleleft S)$ is empty, then $(T - x) \lhd S = (T \lhd S) - x$.  Otherwise, $(T - x) \lhd S = (T \lhd S) - x \sqcup x'$, where $x' = \max(T_{<x}\setminus (T\triangleleft S))$.
\end{lemma}
\begin{proof}
    Let $i$ be such that $T(i) = x$.
    Find smallest 
    $i'$ such that for any $i' \leq j \leq i$,
    $T(j) \in T \lhd S$.
    Notice that $i' > 1$ if and only if
    $T_{<x}\setminus (T\triangleleft S)$ is not empty.
    In this case, $T(i' - 1) = x'$.
    
    Consider the non-recursive way to compute
    $T \lhd S$ and $(T-x) \lhd S$. 
    If $x = T(i)$ is picked by some $s \in S$
    when computing $T \lhd S$,
    then $s$ would pick $T(i-1)$ instead
    when computing $(T-x) \lhd S$.
    Inductively, for $j = i, i-1, \cdots, i' + 1$,
    if $s \in S$
    picks $T(j)$  when computing $T \lhd S$, 
    $s$ would pick $T(j-1)$ instead when computing $(T-x) \lhd S$.
    If $i' = 1$,
    we have $(T-x) \lhd S = T \lhd S - x$
    and we are done. 
    Otherwise, when $s \in S$
    picks $T(i')$ when computing
    $T \tl S$,
    $s$ would pick $T(i'-1) = x'$ instead when computing $(T-x) \lhd S$.
    After that,
    the computations of $T \lhd S$ and $(T-x) \lhd S$
    behave the same.
\end{proof}

\subsection{Domination}

Clearly, $|T \lhd S| \leq |S|$. 
In Example~\ref{E: lhd},
we see sometimes $|T \lhd S| < |S|$.
We would like to understand when $T \lhd S$ has 
the same size as $S$. 
This motivation brings us to the following notion:
\begin{defn}
If $|T|\geq |S|$ and $T(i)<S(i)$ for all $i\in [|S|]$, we say that $S$ \definition{dominates} $T$ and write $T\preceq S$.     
\end{defn}

An alternate characterization for dominance is that $S_1\preceq S_2 \preceq \dots \preceq S_+$ if and only if there exists an $n$-column increasing tableau of normal shape with column $i$ consisting of the numbers in $S_i$.
From this perspective, it is clear that if $T \preceq S$ and $S' \subseteq S$, then $T \preceq S'$.
Notice that $\preceq$ is transitive, but it is not irreflexive since
$\emptyset \preceq \emptyset$.

This notion characterizes when $T \triangleleft S$
has the same size as $S$:
\begin{lemma}
\label{lemma:dominate}
We have $|T\triangleleft S| = |S|$ if and only if $T \preceq S$.
\end{lemma}
\begin{proof}
Prove by induction on $|S|$.  If $|S| = 0$, then $S= \emptyset$ and 
$T\triangleleft \emptyset = \emptyset$.
Correspondingly, we have $T \preceq \emptyset$.

Now assume that the statement is true when $|S| = k-1$, 
and say that $|S| = k$ for some $k \geq 1$.  If $\max(S) \leq \min(T)$, 
then $T\triangleleft S = \emptyset$
and we cannot have $T \preceq S$,
completing the proof in this case.  
Otherwise, $\max(S) > \min(T)$. 
Let $m = \max(T_{< \max(S)})$,
so $T\triangleleft S = m \sqcup (T_{< \max(S)} \triangleleft (S-\max(S)))$.
By the inductive hypothesis, 
we deduce:
\begin{align*}
    |T\triangleleft S| = k &\Leftrightarrow  |T_{< \max(S)} \triangleleft (S-\max(S))|  = k-1\\
    & \Leftrightarrow T_{< \max(S)} \preceq (S-\max(S))\\
    & \Leftrightarrow |T_{< \max(S)}| \geq k-1,\textrm{ and for each $i \in [k-1]$, } T_{< \max(S)}(i) \leq S(i)\\
    & \Leftrightarrow |T| \geq k,\textrm{ and for each $i \in [k]$, } T(i) \leq S(i)\\
    & \Leftrightarrow T\preceq S.\qedhere
\end{align*}
\end{proof}

We are mainly interested in $T \triangleleft S$ 
when $T \preceq S$.
In this case, $|T \lhd S| = |S|$.
If we consider the non-recursive way
to compute $T \lhd S$,
when $S(i)$ picks a number in $T$,
we know it is $(T \lhd S)(i)$.
Thus, we can use the non-recursive
definition to study $T \lhd S$.
We first show that $T \lhd S$
has an entry-wise upper bound. 

\begin{cor}
\label{cor:dominate2}
If $T\preceq S$, then $T\triangleleft S\preceq S$.
\end{cor}
\begin{proof}
Since $T\preceq S$, $|S|=|T\triangleleft S|$, so $|S|\leq|T\triangleleft S|$.
Using the non-recursive definition of $\triangleleft$, the number $S(i)$ picks someone that is smaller than $S(i)$, so $(T\triangleleft S)(i)<S(i)$, completing the proof.
\end{proof}

Using Corollary~\ref{cor:dominate2},
we deduce a simple result
that allows us to understand 
the small numbers in $T \tl S$ in certain cases.

\begin{lemma}
\label{lemma:horizontalmap1}
Assume $T\preceq S$.
Take $x$ such that $|T_{<x}| = |S_{\leq x}|$.
Then for any $j \leq |T_{<x}|$, we have
$(T\triangleleft S)(j) = T(j)$.
\end{lemma}
\begin{proof}
Let $i = |T_{<x}|$.
If $i = 0$, our statement is trivial.
Now assume $i > 0$.
By $T \preceq S$ and Corollary~\ref{cor:dominate2},
$(T \tl S)(i) < S(i) \leq x$,
so there are at least $i$ numbers in $T \tl S$
that are less than $x$.
There are exactly $i$ such numbers in $T$,
so they are all in $T \tl S$.
In other words, for any $j \in [i]$,
$(T\triangleleft S)(j) = T(j)$.
\end{proof}

\begin{cor}
\label{lemma:insertline}
Assume $T\preceq S$ and there exists $x\not\in S$ such that $|S_{<x}| = |T_{<x}|$.
Take $y\in S$ such that $x> y$
and let $S' = (x\sqcup S) - y$.
Then $T\preceq S'$ and for all $i$ such that $S(i)< x$, $(T\triangleleft S)(i) = (T\triangleleft S')(i) = T(i)$. 
\end{cor}
\begin{proof}
For each $i \in [|S|]$,
we have $T(i) < S(i) \leq S'(i)$.
Thus, $T \preceq S'$.
Notice that $|S_{\leq x}| = |S'_{\leq x}| = |T_{<x}|$.
Then the lemma is finished by
Lemma~\ref{lemma:horizontalmap1}.
\end{proof}

Using results from the previous section,
we can get another lower bound 
of $T\lhd S$ involving $(T -x) \lhd S$.

\begin{lemma}
\label{lemma:singlereduction2}
Assume $T\preceq S$.  Take any $x\in T$ such that $T-x\preceq S$.
Then $(T\triangleleft S)(i) \geq ((T-x)\triangleleft S)(i)$ for $i \in [|S|]$.
\end{lemma}

\begin{proof}
    If $x\not\in T\triangleleft S$, then by Lemma \ref{lemma:singlereduction} $T\triangleleft S = (T-x) \triangleleft S$, and the statement is true.
    Thus, we assume $x\in T\triangleleft S$.
    
    Let $x' \in T$ be the largest such that
    $x' < x$ and $x' \notin T \tl S$.
    It $x'$ does not exist,
    by Lemma \ref{lemma:singlereduction3},
    $$|(T-x) \triangleleft S| = |(T\triangleleft S)-x| < |S|.$$
    By Lemma~\ref{lemma:dominate} and $(T-x) \preceq S$,
    we cannot have $|(T-x) \preceq S| < |S|$. 
    Thus, $x'$ must exist.
    By Lemma \ref{lemma:singlereduction3} $(T - x) \lhd S = (T \lhd S) - x \sqcup x'$.
    Our statement is proved since $x' < x$.
\end{proof}

We end this section with a proposition that allows us to ``split'' the  
computation of $T \tl S$.

\begin{prop}
\label{P: Split}
Find $x$ such that $T_{\geq x} \preceq S_{>x}$.
Then 
$$
T \tl S = (T_{< x} \tl S_{ \leq x}) 
\sqcup (T_{\geq x} \tl S_{>x})
$$
\end{prop}
\begin{proof}
Consider the non-recursive way
of computing $T \tl S$
and $(T_{\geq x} \tl S_{>x})$.
Since $T_{\geq x} \preceq S_{>x}$,
$s \in S_{>x}$ will pick the same 
$t \in T_{\geq x}$ in both computations.

Now consider the computation of 
$T \tl S$.
After every $s \in S_{>x}$ picks,
the non-picked numbers in $T$
are $T_{<x}$ together with
some numbers at least $x$.
When the next $s \in S$ picks, 
we know $s \leq x$, 
so it will only pick numbers from $T_{<x}$.
Starting from here,
the computation behaves the same
as if computing $T_{<x} \tl S_{\leq x}$.
\end{proof}

\begin{cor}
\label{C: n bullet} 
Find $x \in T$ such that $T_{\geq x+1} \preceq S_{>x + 1}$ and $x + 1 \notin S$.
Then $x \notin T \tl S$.
\end{cor}
\begin{proof}
By Proposition~\ref{P: Split},
$T \tl S = (T_{< x + 1} \tl S_{ \leq x + 1}) 
\sqcup (T_{\geq x + 1} \tl S_{>x + 1})$.
Clearly $x$ is not in the second term. 
By $x+1 \notin S$,
$x$ is not in the first term either.
\end{proof}

\subsection{Comparing $T \tl S'$
and $T \tl S$ when $S' \subseteq S$}

We start with a simple property.

\begin{lemma}
\label{lemma:subset1}
For $S'\subseteq S$, 
we have
$T\triangleleft S' \subseteq T \triangleleft S$.
\end{lemma}
\begin{proof}
Prove by induction on $|S|$.  If $|S|=0$, then $S = S' = \emptyset$, so $T\triangleleft S' = \emptyset = T\triangleleft S$.

Now assume that the statement is true whenever $|S|=k-1$, and consider an arbitrary $S$ such that $|S|=k$.  If $\max(S) \leq \min(T)$, then  $T\triangleleft S' = T\triangleleft S = \emptyset$ and we are done.  
Otherwise, let $m = \max(T_{<\max(S)})$.
If $m \in S'$, by the inductive hypothesis, 
$$
T \tl S' 
= m \sqcup (T_{< m} \tl (S'-m))
\subseteq m \sqcup (T_{< m} \tl (S-m))
= T \tl S.
$$
Otherwise, $m \notin S'$.
We have $S' \subseteq (S - m)$.
Notice that $T \tl S \subseteq T_{<m}$.
By Corollary \ref{lemma:Reduction (subset)}
and the inductive hypothesis,
$$
T \tl S' 
= T_{< m} \tl S' 
\subseteq (T_{< m} \tl (S-m))
= T \tl S.
$$
\end{proof}

By Lemma \ref{lemma:subset1}, when $T\preceq S$, $T\triangleleft (S-a)$ will be $T\triangleleft S$ with one element removed. For certain $a$,
we can figure out which element of $T \tl S $ is removed.

\begin{lemma}
\label{lemma:subsetkeycalc}
Assume $T\preceq S$ and take $S'\subset S$. If $a = \min(S\setminus S')$ and $b = \min((T\triangleleft S) \setminus (T\triangleleft S'))$, then $T\triangleleft (S-a)= (T\triangleleft S)-b$
\end{lemma}

\begin{proof}
Since $S-a$ and $S'$ are subsets
of $S$,
by Lemma~\ref{lemma:subset1},
$T \tl (S-a)$ and $T \tl S'$
are subsets of $T \tl S$.
By Corollary~\ref{lemma:Reduction (subset)},
we may replace $T$ by $T \tl S$.
This change does not affect $T \tl S$,
$T \tl (S-a)$ and $T \tl S'$.

By Lemma \ref{lemma:dominate}, $|T\triangleleft S| = |S| = 1+ |S-a| = 1 + |T\triangleleft (S-a)|$. By Lemma \ref{lemma:subset1} we have that $T\triangleleft (S-a) \subseteq T\triangleleft S$.
Thus, $(T\triangleleft S) \setminus (T\triangleleft (S-a))$ contains only one element, which we denote as $x$. It suffices to show that $x = b$.
Let $x = T(i)$.
It remains to check: $x \notin T \tl S'$
and $T(j) \in T \tl S'$ for any $j \in [i-1]$.
\begin{itemize}
\item Since $x \notin T \tl (S-a)$
and $S' \subseteq S-a$,
by Lemma~\ref{lemma:subset1},
$x \notin T \tl S'$.
\item On one hand, we have 
$S(i) > T(i) =x$.
On the other hand, 
consider the non-recursive way
to compute $T \tl (S - a)$.
When $(S-a)(i-1)$ picks $T(i-1)$,
$x$ has not been picked. 
Thus, $(S-a)(i-1) \leq x$.
In particular, 
$S(i) > (S-a)(i-1)$,
so $S(i) \leq a$
and $(S-a)(i-1) = S(i-1)$.

Now we have $|S_{\leq x}| = i-1$.
By the definition of $a$,
$S(i) \leq a$ implies
that $S(1), \cdots, S(i-1) \in S'$,
so $|S'_{\leq x}| = i-1 = |T_{<x}|$.
By Lemma~\ref{lemma:horizontalmap1},
$(T \tl S')(j) = T(j)$ for $j \in [i-1]$.
\end{itemize}
\end{proof}

\subsection{Understanding $U \tl T \tl S$}
The main goal of this section is to prove the following technical lemma:
\begin{lemma}
\label{lemma:insert7}
Suppose $U\preceq T\preceq S$.  Assume that there exists $x\not\in T$ such that $|U_{<x}| = |T_{<x}|$.  Assume there exists $y \in T_{<x}$ such that $T \sqcup x - y \preceq S$. Take the smallest such $y$, then $$U \tl T \tl S' = U \tl (T \sqcup x - y) \tl S'$$ where $S'$ is any subset of $S$.
Consequently, 
$$
U \tl T = U \tl (T \sqcup x - y).
$$
\end{lemma}

To prove this Lemma, 
we consider three sets:
\begin{align}
\label{EQ: Three sets}
U \lhd T \lhd S',\quad
U \lhd (T\sqcup x) \lhd S', \quad
U \lhd (T \sqcup x - y) \lhd S'.
\end{align}

The goal is to show the first set equals 
the last set in~(\ref{EQ: Three sets}.
We first derive a lemma that implies 
the last two agree.

\begin{lemma}
\label{lemma:deletion7}
Assume $T\preceq S$ and $|T|>|S|$. Let $y \in T$ be the smallest
such that $(T-y) \preceq S$. 
Then 
\begin{itemize}
\item For all $T(j)<y$, $(T\triangleleft S)(j) = T(j)$.
\item We have $y \notin T \lhd S$.
\item For any $S'\subseteq S$, $(T-y)\triangleleft S' = T\triangleleft S'$.
\end{itemize}

\end{lemma}
\begin{proof}
Let $i$ be such that $y = T(i)$. 
\begin{itemize}
\item If $i=1$, then $y = \min(T)$, and 
the first statement is vacuously true. Now assume that $i>1$. By the definition of $y$, $T-T(i)\preceq S$, but $T-T(i-1)\not\preceq  S$. 
Thus, $T(i)\geq S(i-1)$.  
In other words, $|S_{\leq y}| = i-1 = |T_{<y}|$.
Then the first statement follows from Lemma~\ref{lemma:horizontalmap1}.
\item By the first statement,
for all $j< i$, $(T\triangleleft S)(j) = T(j) < T(i) = y$.  For all $j\geq i$, Lemma \ref{lemma:singlereduction2} says that $(T\triangleleft S)(j) \geq ((T-y)\triangleleft S)(j) \geq (T-y)(j) = T(j+1) > T(i) = y$.  Therefore, for all $j$, $(T\triangleleft S)(j)\not= y$, so $y\not\in T\triangleleft S$.
\item By Lemma \ref{lemma:subset1}, $T\triangleleft S' \subseteq T\triangleleft S$. By the second statement, 
$y\not\in T\triangleleft S$, so $y\not\in T\triangleleft S'$. 
By Lemma \ref{lemma:singlereduction}, $(T-y)\triangleleft S' = T\triangleleft S'$.
\end{itemize}
\end{proof}

To show the first two sets in~(\ref{EQ: Three sets}) agree, we need the following two lemmas:

\begin{lemma}
\label{L: reverse complement}
Suppose $|U| \leq |T|$.
Let $\delta = |U| - |T|$.
Suppose for each $i \in [|U|]$,
$U(i) \leq T(i + \delta)$.
Take $x$ such that $x > \max(T)$.
\begin{itemize}
\item We have $U \tl (T \sqcup x) = U \tl T = T$.
\item For any $T' \subsetneq T$, let $x' = \max(T \setminus T')$.
Then $U \tl (T' \sqcup x) = U \tl (T' \sqcup x')$.
\end{itemize}
\end{lemma}
\begin{proof}
The first statement is immediate.
For the second statement,
assume $x' = T(i)$.
We know $T(|T|) ,\cdots,  T(i+1)\in T'$.
Consider the non-recursive way
to compute $U \tl (T' \sqcup x)$ 
and $U \tl (T' \sqcup x')$.
In the first process, clearly, 
$x, T(|T|) ,\cdots,  T(i+1)$
would pick $U(|U|), U(|U| - 1) \cdots, U(i - \delta)$. 
In the second process,
$T(|T|),\cdots, T(i)$
would also pick 
$U(|U|), \cdots, U(i - \delta)$.
After that, 
the two processes have the same behavior. 
\end{proof}

\begin{lemma}
\label{L: UT}
Assume we can find $x \notin T$
such that $U_{\geq x} \preceq T_{> x}$ and 
$|U_{<x}| \leq |T_{< x}|$.
Let $\delta = |T_{< x}| - |U_{<x}|$.
Assume for any $i \in [|U_{<x}|]$,
$U(i) \leq T(i + \delta)$.
Then $U \lhd T \lhd S = U \lhd (T \sqcup x) \lhd S$ for any set $S$.
\end{lemma}
\begin{proof}
Suppose $x \notin (T \sqcup x) \lhd S$. 
By Lemma~\ref{lemma:singlereduction},
$(T \sqcup x) \lhd S = T \lhd S$
so our claim is trivial.

Otherwise, 
let $T' = (T \sqcup x) \lhd S$. 
Since $U_{\geq x} \preceq T_{> x}$,
we have $U_{\geq x} \preceq T'_{>x}$.
By Proposition~\ref{P: Split},
\begin{align}
\label{EQ: UT}
U \tl T' = U_{<x} \tl T'_{\leq x} \sqcup U_{\geq x} \preceq T'_{> x}.    
\end{align}

By Lemma~\ref{lemma:singlereduction3}, $(T \tl S)_{>x} = T'_{>x}$.
Thus, 
~(\ref{EQ: UT}) still holds
if we replace $T'$ by $T \tl S$.
To show $U \tl T' = U \tl (T \tl S)$,
it remains to check:
\begin{align}
\label{EQ: UT2}
U_{<x} \tl T'_{\leq x} = U_{<x} \tl 
(T \tl S)_{\leq x}.    
\end{align}

Let $x' = \max(T_{<x} \setminus T'_{<x})$.
If $x'$ does not exist,
by Lemma~\ref{lemma:singlereduction3},
$T'_{<x}= (T \tl S)_{<x} \sqcup x$.
Then~(\ref{EQ: UT2}) follows from
the first statement of 
Lemma~\ref{L: reverse complement}.
Otherwise, 
by Lemma~\ref{lemma:singlereduction3},
$T'_{<x}= (T \tl S)_{<x} \sqcup x - x'$.
Then~(\ref{EQ: UT2}) follows from
the second statement of 
Lemma~\ref{L: reverse complement}.
\end{proof}

Finally, we can prove Lemma~\ref{lemma:insert7}
\begin{proof}[Proof of Lemma~\ref{lemma:insert7}]
By the third statement of Lemma~\ref{lemma:deletion7},
the last two sets in~(\ref{EQ: Three sets}) agree.
We just need to check the first two sets 
in~(\ref{EQ: Three sets}) agree. 
Since $U \preceq T$ and $|U_{<x}| = |T_{<x}|$,
we have $U_{\geq x} \preceq T_{\geq x}$.
Then $U$ and $T$ satisfies the condition
in Lemma~\ref{L: UT},
so we have the first equation in Lemma~\ref{lemma:insert7}.
For the second equation, 
we simply let $S$ be a set consisting
of $|T|$ very large numbers,
so $T \lhd S = T$ 
and $(T \sqcup x - y) \lhd S = (T \sqcup x - y)$.\qedhere
\end{proof}

\section{The left key of an increasing tableau}
\label{S: left key}
For an increasing tableau $P$,
its left key is denoted as $K_-(P)$.
In Section~\ref{SS: Kjdt}, 
we give the usual definition of $K_-(P)$ as in~\cite{RY}
using K-theoretic jeu-de-taquin of Thomas and Yong~\cite{TY}. 
In Section\ref{SS: left key and tl},
we derive a simple way to compute $K_-(P)$
using the $\tl$ operator. 

\subsection{Defining $K_-(P)$ using K-theoretic jeu-de-taquin}
\label{SS: Kjdt}

In this section, 
we consider a more general family of tableaux.
Let $\lambda$ and $\mu$ be two partitions
such that the Young diagram of $\mu$ is contained in the
Young diagram of $\lambda$.
The \definition{skew Young diagram} of 
$\lambda/\mu$ is the set theoretic
difference between the Young diagram of $\lambda$
and the Young diagram of $\mu$.
A \definition{skew tableau} of shape $\lambda/\mu$ 
is a filling of the skew
Young diagram of $\lambda/\mu$
with positive integers.
Skew tableaux whose shapes are 
left-justified and top-justified 
are called \definition{normal}.
Skew tableaux whose cells are right-justified 
and bottom-justified 
are called \definition{anti-normal}.
We say a skew tableau is increasing
if each of its row (resp. column)
increases from left to right (resp. top to bottom).
Let $\sInc$ be the set of all increasing tableaux
with skew shape.

Thomas and Yong~\cite{TY} introduce an operation
called \definition{anti-rectification}
which shifts cells in $T \in \sInc$ to
make it anti-normal. 
This map consists of small moves called
\definition{reverse K-jeu-de-taquin} (\definition{$\revKjdt$})
which are defined on 
a more general family of fillings.
\begin{defn}

A \definition{ dotted skew tableau} is a filling of a skew Young
diagram with positive integers and the symbol~$\bullet$.
Let $m$ be a positive integer. 
Let $\sInc_{m < \bullet}$ be the set of dotted skew tableaux
such that each row and column is strictly increasing
with respect to the order 
$\cdots < m < \bullet < m+1 < \cdots$.
\end{defn}

The shape of a dotted skew tableau is defined similarly.
Therefore, we can also use normal and anti-normal
to describe dotted skew tableaux.

\begin{defn}~\cite{TY}
Define the \definition{reverse K-jeu-de-taquin} 
(\definition{$\revKjdt$})
as a shape-preserving map 
from $\sInc_{m < \bullet}$ 
to $\sInc_{m - 1 < \bullet}$
when $m \geq 1$.
Take $T \in \sInc_{m < \bullet}$.
For each $\bullet$ (resp. $m$)
that is adjacent to an $m$ (resp. $\bullet$) in $T$, 
we replace it by $m$ (resp. $\bullet$).
The resulting dotted skew tableau is in 
$\sInc_{m-1 < \bullet}$.
\end{defn}

\begin{exa}
We perform $\revKjdt$ to send an element of 
$\sInc_{4 < \bullet }$ to $\sInc_{3 < \bullet}$
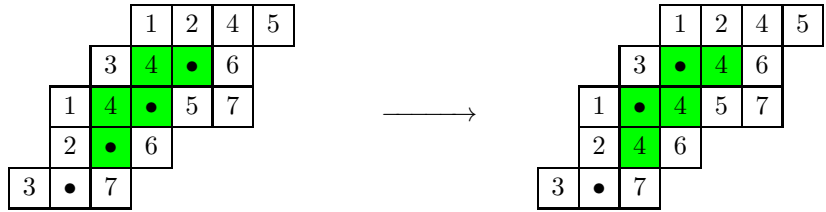
$$
\raisebox{1cm}{\begin{ytableau}
\none & \none & \none   & 1   & 2   & 4 & 5 \\
\none & \none & 3   & *(green)4      & *(green)\bullet & 6 \\
\none & 1     & *(green)4   &*(green)\bullet & 5 & 7 \\
\none & 2     & *(green)\bullet   &6 \\
3 & \bullet     & 7\\
\end{ytableau}}
\quad \quad \quad \xrightarrow{\quad\quad\quad} \quad \quad 
\raisebox{1cm}{\begin{ytableau}
\none & \none & \none   & 1   & 2   & 4 & 5 \\
\none & \none & 3   & *(green)\bullet      & *(green)4 & 6 \\
\none & 1     & *(green)\bullet   &*(green)4 & 5 & 7 \\
\none & 2     & *(green)4   &6 \\
3 & \bullet     & 7\\
\end{ytableau}}
$$
\end{exa}

We provide another way to 
understand $\revKjdt$.
Let us start with following definitions.

\begin{defn}
A skew shape is called a {\bf short ribbon}
if it satisfies the following.
\begin{itemize}
\item All cells are connected. 
\item Does not have a $2 \times 2$ subshape.
\item Each column or row has at most 2 cells. 
\end{itemize}

A dotted skew tableau is called
an {\bf alternating ribbon} of $m$
if it satisfies the following.
\begin{itemize}
\item Its shape is a short ribbon. 
\item All cells are filled by $\bullet$ and $m$.
Adjacent cells are filled differently. 
\end{itemize}
\end{defn}

It is clear that for any short ribbon, 
there are two alternating ribbons of $m$
with that shape.
We say an alternating ribbon has type 1
if its topmost cell 
in the leftmost column is $m$.
Otherwise, it has type 2. 

\begin{exa}
In this example, we present a short ribbon
and the two alternating ribbons of $m$
with that shape.
The former has type 1 and the latter has type 2. 

$$
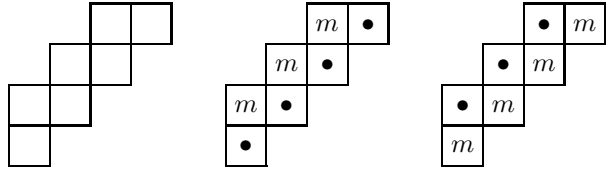
\begin{ytableau}
\none & \none &   &  \\
\none &   &   \\
  &  \\
 & \none
\end{ytableau}
\quad \quad
\begin{ytableau}
\none & \none & m& \bullet\\
\none & m & \bullet \\
m & \bullet\\
\bullet
\end{ytableau}
\quad \quad
\begin{ytableau}
\none & \none & \bullet & m\\
\none & \bullet & m \\
\bullet & m\\
m
\end{ytableau}
$$
\end{exa}

Now take $T \in \sInc_{m < \bullet}$.
We may describe $\revKjdt$ as follows.
If we focus on its $m$ and~$\bullet$,
ignoring all other numbers,
we see a skew-shape tableau
where each connected component is an
alternating ribbon.
Each alternating ribbon
with more than one cell
must have type 1. 
A $\revKjdt$ move would change these 
alternating ribbons
with more than one cell into type 2. 

Next, we use $\revKjdt$ to describe a process called \definition{anti-rectification}. 
Let $T$ be an increasing tableau 
of skew shape $\lambda/\mu$.
Say all numbers in $T$ are at most $m$ for some $m$.
We may add some cells containing $\bullet$ to $T$
and obtain an element of $\sInc_{m < \bullet}$.
By recursively applying $\revKjdt$, we obtain 
an element of $\sInc_{0 < \bullet}$.
Then we remove all $\bullet$.
This process yields a skew increasing tableau
with a different shape. 
By repeating this process,
we will end with an increasing tableau 
of anti-normal shape. 

\begin{exa}
\label{E: anti-rectify}
Let $T$ be the following element of $\sInc$.
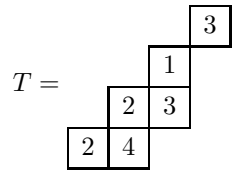
$$
T = \raisebox{0.6cm}{\begin{ytableau}
\none & \none & \none  &  3\\
\none & \none & 1\\
\none & 2 & 3\\
2 & 4\\
\end{ytableau}}
$$
We may put a $\bullet$ under each "3" and perform
$\revKjdt$:
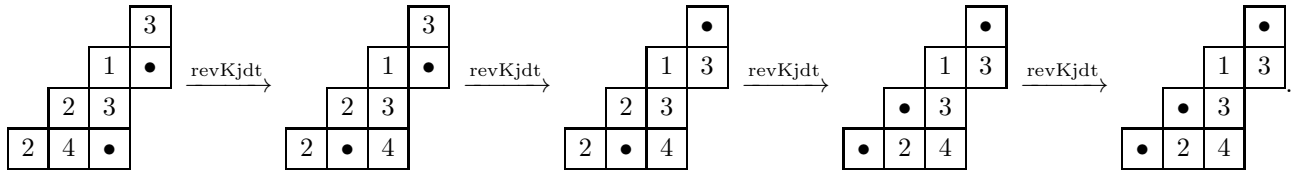
$$
\raisebox{0.6cm}{\begin{ytableau}
\none & \none & \none  &  3\\
\none & \none & 1 & \bullet\\
\none & 2 & 3\\
2 & 4 & \bullet\\
\end{ytableau}}
\: \xrightarrow{\revKjdt} \:
\raisebox{0.6cm}{\begin{ytableau}
\none & \none & \none  &  3\\
\none & \none & 1 & \bullet\\
\none & 2 & 3\\
2 & \bullet & 4\\
\end{ytableau}}
\: \xrightarrow{\revKjdt} \:
\raisebox{0.6cm}{\begin{ytableau}
\none & \none & \none  &  \bullet\\
\none & \none & 1 & 3\\
\none & 2 & 3\\
2 & \bullet & 4\\
\end{ytableau}}
\: \xrightarrow{\revKjdt} \:
\raisebox{0.6cm}{\begin{ytableau}
\none & \none & \none  &  \bullet\\
\none & \none & 1 & 3\\
\none & \bullet & 3\\
\bullet & 2 & 4\\
\end{ytableau}}
\: \xrightarrow{\revKjdt} \:
\raisebox{0.6cm}{\begin{ytableau}
\none & \none & \none  &  \bullet\\
\none & \none & 1 & 3\\
\none & \bullet & 3\\
\bullet & 2 & 4\\
\end{ytableau}}.
$$
Now we may ignore the $\bullet$,
obtaining an element of $\sInc$.
Then we may place a $\bullet$
at row $3$ and repeat this process:
$$
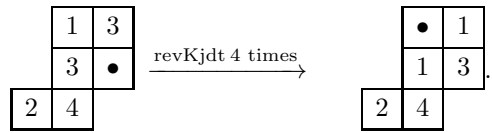
\raisebox{0.4cm}{\begin{ytableau}
\none & \none & 1 & 3\\
\none & \none & 3 & \bullet\\
\none & 2 & 4\\
\end{ytableau}}
\: \xrightarrow{\revKjdt \:4 \textrm{ times}} \:
\raisebox{0.4cm}{\begin{ytableau}
\none & \none & \bullet & 1\\
\none & \none & 1 & 3\\
\none & 2 & 4\\
\end{ytableau}}.
$$

Finally, we remove the current $\bullet$ and 
put a $\bullet$ under 3.
We perfrom $\revKjdt$ again, 
obtaining
$$
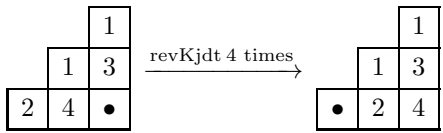
\raisebox{0.4cm}{\begin{ytableau}
\none & \none & \none & 1\\
\none & \none & 1 & 3\\
\none & 2 & 4 & \bullet\\
\end{ytableau}}
\: \xrightarrow{\revKjdt \:4 \textrm{ times}} \:
\raisebox{0.4cm}{\begin{ytableau}
\none & \none & 1\\
\none & 1 & 3\\
\bullet & 2 & 4\\
\end{ytableau}}
$$
which is the result of our anti-rectification
if we ignore its $\bullet$
\end{exa}

\begin{rem}
During the anti-rectification process, 
we need to choose where to place the $\bullet$.
Making different choices will actually affect the final
result. 
Let us anti-rectify the $T$ in Example~\ref{E: anti-rectify}
with different choices of the positions to place the $\bullet$.
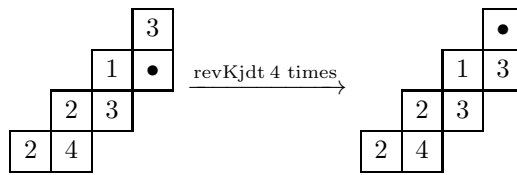
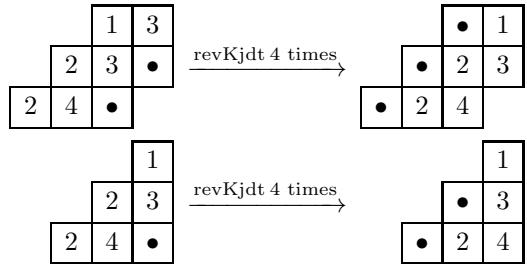
$$
\raisebox{0.6cm}{\begin{ytableau}
\none & \none & \none  &  3\\
\none & \none & 1 & \bullet\\
\none & 2 & 3\\
2 & 4 \\
\end{ytableau}}
\: \xrightarrow{\revKjdt \:4 \textrm{ times}} \:
\raisebox{0.6cm}{\begin{ytableau}
\none & \none & \none  &  \bullet\\
\none & \none & 1 & 3\\
\none & 2 & 3\\
2 & 4 \\
\end{ytableau}}
$$
$$
\raisebox{0.4cm}{\begin{ytableau}
\none & \none & 1 & 3\\
\none & 2 & 3 & \bullet\\
2 & 4 & \bullet \\
\end{ytableau}}
\: \xrightarrow{\revKjdt \:4 \textrm{ times}} \:
\raisebox{0.4cm}{\begin{ytableau}
\none & \none & \bullet & 1\\
\none & \bullet & 2 & 3\\
\bullet & 2 & 4 \\
\end{ytableau}}
$$
$$
\raisebox{0.4cm}{\begin{ytableau}
\none & \none & \none & 1\\
\none & \none & 2 & 3\\
\none & 2 & 4 & \bullet \\
\end{ytableau}}
\: \xrightarrow{\revKjdt \:4 \textrm{ times}} \:
\raisebox{0.4cm}{\begin{ytableau}
\none & \none & \none & 1\\
\none & \none & \bullet & 3\\
\none & \bullet & 2 & 4 \\
\end{ytableau}}
$$
If we ignore the bullet in the last dotted skew tableau,
we obtain the result of our anti-rectification. 
This time, we end up with something different
from the output in Example~\ref{E: anti-rectify}.
\end{rem}

Finally, for an increasing tableau $P$
of normal shape.
We use anti-rectification to define 
$K_-(P)$, the left key of $P$. 
Let $P_{\leftarrow j}$
to be the tableau obtained by keeping
only the first $j$ columns of $P$.
Column $j$ of $K_-(P)$ only depends on $P_{\leftarrow j}$.
Say $P$ has $r$ rows.
We may embed $P_{\leftarrow j}$ in an $j$ by
$r$ rectangle.
Then we anti-rectify $P_{\leftarrow j}$.
When we choose where to place the $\bullet$,
we always place one $\bullet$ at the leftmost available 
position inside the rectangle. 
After the anti-rectification,
numbers in the leftmost column 
will form column $j$ of $K_-(P)$.

\begin{exa}
Let $P$ be the following increasing tableau
of normal shape. 
$$
P = \raisebox{0.6cm}
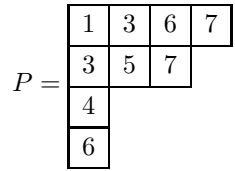
{\begin{ytableau}
1 & 3 & 6  &  7\\
3 & 5 & 7\\
4\\
6\\
\end{ytableau}}
$$
We compute $K_-(P)$ as follows. 
First, consider $P_{\leftarrow 1}$,
which is already anti-normal. 
Column 1 of $K_-(P)$ consists of $\{1, 3, 4, 6\}$.
Next, consider $P_{\leftarrow 2}$.
We may anti-rectify it as follows.
$$
P_{\leftarrow 2} = \raisebox{0.6cm}
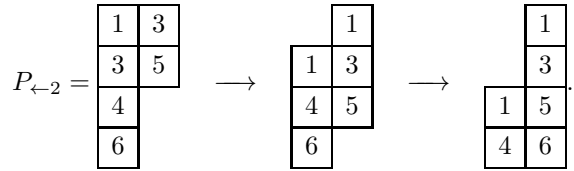
{\begin{ytableau}
1 & 3\\
3 & 5\\
4\\
6\\
\end{ytableau}}
\quad \xrightarrow{\quad} \quad 
\raisebox{0.6cm}{\begin{ytableau}
\none & 1\\
1 & 3\\
4 & 5\\
6\\
\end{ytableau}}
\quad \xrightarrow{\quad} \quad 
\raisebox{0.6cm}{\begin{ytableau}
\none & 1\\
\none & 3\\
1 & 5\\
4 & 6\\
\end{ytableau}}.
$$
Column $2$ of $K_-(P)$ consists
of $\{1,4\}$.
Next, we anti-rectify $P_{\leftarrow 3}$.
Notice that the first two iterations
of anti-rectifying $P_{\leftarrow 3}$
will be the same as anti-rectifying $P_{\leftarrow 2}$.
$$
P_{\leftarrow 3} = \raisebox{0.6cm}
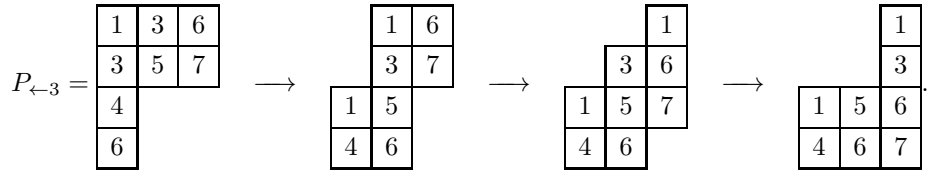
{\begin{ytableau}
1 & 3 & 6\\
3 & 5 & 7\\
4\\
6\\
\end{ytableau}}
\quad \xrightarrow{\quad} \quad
\raisebox{0.6cm}{\begin{ytableau}
\none & 1 & 6\\
\none & 3 & 7\\
1 & 5\\
4 & 6\\
\end{ytableau}}
\quad \xrightarrow{\quad} \quad
\raisebox{0.6cm}{\begin{ytableau}
\none & \none & 1\\
\none & 3 & 6\\
1 & 5 & 7\\
4 & 6\\
\end{ytableau}}
\quad \xrightarrow{\quad} \quad
\raisebox{0.6cm}{\begin{ytableau}
\none & \none & 1\\
\none & \none & 3\\
1 & 5 & 6\\
4 & 6 & 7\\
\end{ytableau}}.
$$

Column $3$ of $K_-(P)$ also consists
of $\{1, 4\}$.
Finally, we anti-rectify $P_{\leftarrow 4}$.
$$
P_{\leftarrow 4} = \raisebox{0.6cm}
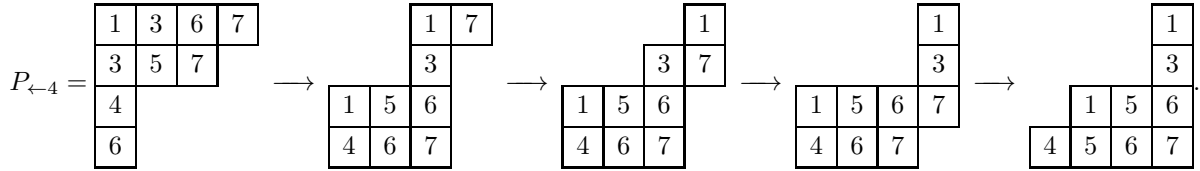
{\begin{ytableau}
1 & 3 & 6 & 7\\
3 & 5 & 7\\
4\\
6\\
\end{ytableau}}
\: \xrightarrow{\quad} \:
\raisebox{0.6cm}{\begin{ytableau}
\none & \none & 1 & 7\\
\none & \none & 3\\
1 & 5 & 6\\
4 & 6 & 7\\
\end{ytableau}}
\: \xrightarrow{\quad} \:
\raisebox{0.6cm}{\begin{ytableau}
\none & \none & \none & 1\\
\none & \none & 3 & 7\\
1 & 5 & 6\\
4 & 6 & 7\\
\end{ytableau}}
\: \xrightarrow{\quad} \:
\raisebox{0.6cm}{\begin{ytableau}
\none & \none & \none & 1\\
\none & \none & \none & 3\\
1 & 5 & 6 & 7\\
4 & 6 & 7\\
\end{ytableau}}
\: \xrightarrow{\quad} \:
\raisebox{0.6cm}{\begin{ytableau}
\none & \none & \none & 1\\
\none & \none & \none & 3\\
\none & 1 & 5 & 6\\
4 & 5 & 6 & 7\\
\end{ytableau}}.
$$
Column $4$ of $K_-(P)$ consists
of $\{4\}$.
Thus, $K_-(P)$ is 
$$
P = \raisebox{0.6cm}
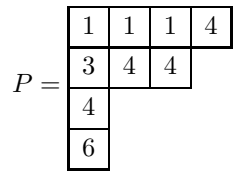
{\begin{ytableau}
1 & 1 & 1  &  4\\
3 & 4 & 4\\
4\\
6\\
\end{ytableau}}
$$
\end{exa} 

From this definition, 
it is not clear that $K_-(P)$
must be a key. 
A proof of this can be found
in~\cite{RY}.
This argument is originally
formed by the first author
of this paper.

It is hard to compute
$K_-(P)$ via this definition
since the computation consists of 
many $\revKjdt$ moves.
Thus, 
in the next section, 
we provide a simple way to compute $K_-(P)$
using $\tl$.

\subsection{Computing $K_-(P)$ using $\tl$}
\label{SS: left key and tl}
We now present a simple method that computes 
$K_-(P)$ for $P \in \Inc_n$.
Recall that $P_i$ is the set of numbers that
appear in column $i$ of $P$.
Define
$$P_{i,j} := P_i \lhd P_{i+1} \lhd \cdots \lhd P_j.$$
Then we can state the main result of this section. 
\begin{thm}
\label{thm:anti and lhd}
Let $P \in \Inc_n$ with $C$ columns.
If we anti-rectify $P$ and get $P^\searrow$ whose rightmost column
is in column $C$,
then $P^\searrow_1 = P_{1, C}$.
\end{thm}

A consequence of this result is the
following simple algorithm 
that computes $K_-(P)$.
\begin{cor}
\label{C: leftkeycalc}
Let $P$ be an increasing tableau with normal shape.
Column $i$ of $K_-(P)$ consists of 
$P_{1, i}$.
\end{cor}

\begin{proof}[Proof of Theorem~\ref{thm:anti and lhd}]
Since $P^\searrow$ is anti-normal,
$P^\searrow_{1,C} = P^\searrow_1$.
We need to show $P_{1,C} = P^\searrow_{1,C}$.
We obtain $P^\searrow$ from
$P$ by applying the 
three operators iteratively:
Adding $\bullet$ to get a tableau in $\sInc_{n<\bullet}$; Performing $\revKjdt$
to get a tableau in $\sInc_{0<\bullet}$;
Removing all $\bullet$.
It remains to check the  Lemma~\ref{lemma:revKjdt preserves triangle} below.   
\end{proof}

For $P \in \sInc_{n < \bullet}$,
we let $P_i$ be the set of numbers
appearing in column $i$ of $P$, 
ignoring the $\bullet$.
Then $P_{1,C}$ is also well-defined.

\begin{lemma}
\label{lemma:revKjdt preserves triangle}
Take $P \in \sInc_{n < \bullet}$.
Let $C$ be the largest such that $P_C \neq \emptyset$.
Perform a $\revKjdt$ move on $P$ 
and obtain $\tilde{P} \in \Inc_{n - 1 < \bullet}$.
Then $ P_{1, C} = \tilde{P}_{1, C}$.
\end{lemma}

This Lemma can be proved by studying the effect of 
changing one alternating ribbon.

\begin{lemma}
\label{lemma:kjdt on one component}
Take $P \in \sInc_{n < \bullet}$ with $C$ columns.
Assume there is one alternating ribbon of $n$ and $\bullet$
within $P$.
Moreover, assume the alternating ribbon has more than one cell
and appears in every column of $P$.

Now switch this alternating ribbon into type 2
of the same shape.
Let $\tilde{P} \in \sInc_{n - 1 < \bullet}$ be the result. 
Then 
$$
P_1 \tl \cdots \tl P_C \tl S = 
\tilde{P}_1 \tl \cdots \tl \tilde{P}_C \tl S
$$
for any set $S$.
\end{lemma}
\begin{proof}
When $C = 1$, the claim is immediate. 
Assume $C > 1$.
We start with following observations. 
\begin{itemize}
\item If $n \notin P_C$, 
then $\tilde{P}_C = P_C \sqcup n$.
Otherwise, $\tilde{P}_C = P_C$.
\item Similarly, if column 1 of $P$ does not have $\bullet$,
then $\tilde{P}_1 = P_1 - n$.
Otherwise, $\tilde{P}_1 = P_1$.
\item For $1 < c < C$,
$\tilde{P}_c = P_c$.
\end{itemize}

First, we claim
$$
P_{C-1} \lhd P_C \lhd S = P_{C-1} \lhd \tilde{P_C} \lhd S.
$$

The claim is immediate if $P_C = \tilde{P_C}$.
Assume $\tilde{P_C} = P_C \sqcup n$.
Notice that setting
$U = P_{C-1}$, $T = P_C$ and $x = n$ would satisfy the conditions of
Lemma~\ref{L: UT}.
Our claim is implied by
Lemma~\ref{L: UT}.

Now, notice that we are done when column 1 of $P$ contains a $\bullet$.
If so, $\tilde{P}_1 = P_1, \dots, \tilde{P}_{C-1} = P_{C-1}$.
Then $P_1 \tl \cdots \tl P_C \tl S = 
\tilde{P}_1 \tl \cdots \tl \tilde{P}_C \tl S$ by the previous claim.

Otherwise, we assume column 1 of $P$ does not have a $\bullet$,
so $\tilde{P}_1 = P_1 - n$.
Notice that setting $T = P_1$, $S = P_2$ and $x = n$ 
would satisfies the conditions of Corollary~\ref{C: n bullet}.
Thus, $n \notin P_1 \tl P2$.
Now by our first claim,
$$
P_1 \tl \cdots \tl P_C \tl S = 
P_1 \tl \cdots \tl P_{C-1} \tl \tilde{P}_C \tl S.
$$
By $n \notin P_1 \tl P2$ and Lemma~\ref{lemma:subset1},
both sides of the equation above do not contain $n$.
Thus, 
$$
P_1 \tl \cdots \tl P_{C-1} \tl \tilde{P}_C \tl S = 
\tilde{P}_1 \tl P_2 \tl \cdots \tl P_{C-1} \tl \tilde{P}_C \tl S.
$$
Then the proof is finished after replacing all $P_i$
on the right 
by $\tilde{P}_i$.
\end{proof}

Now we can prove Lemma~\ref{lemma:revKjdt preserves triangle}.

\begin{proof}
Look at $n$ and~$\bullet$ in $T$.
Assume there are $J$ alternating ribbons
with more than one cell.
Label them
with $1, 2, \dots, J$ arbitrarily. 
Now let $P^0 = P$ 
and
$P^j$ is obtained from $P^{j-1}$
by switching the alternating ribbon $j$
to type 2. 
Then we know $P^J = \tilde{P}$.
It remains to show
$P_{1,C}^{j-1} = P_{1,C}^{j}$.

Assume the alternating ribbon $j$
spans from column $a_j$ to column $b_j$.
Then $P^{j-1}$ and $P$ only differ 
between these two columns.
Thus, column $a_j, \dots, b_j$ of $P^{j-1}$
form a dotted tableau in $\sInc_{n < \bullet}$.
We may apply Lemma~\ref{lemma:kjdt on one component}
and get 
$$
P_{a_j}^{j-1} \tl \dots \tl P_{b_j}^{j-1} \tl S
= P_{a_j}^{j} \tl \dots \tl P_{b_j}^{j} \tl S
$$
for any set $S$.

Finally, set $S = T_{b_j + 1, C}^{j-1}$
if $b_j < C$.
If $b_j = C$, 
set $S = \{1,\dots, N\}$ for some $N$ large enough.
We have $P_{a_j, C}^{j-1} = P_{a_j, C}^{j}$,
which implies $P_{1, C}^{j-1} = P_{1, C}^{j}$
\end{proof}

\section{The reverse row insertion}
\label{S: insertion}

Recall the main tool we used is a bijection
$\Psi : \T \rightarrow \C$ of Huang, Shimozono and Yu.
This map relies on a reverse row insertion algorithm
on increasing tableaux.
In Section~\ref{SS: reverse insertion}, 
we describe the reverse row insertion
algorithm which leads to the bijection $\Psi: \C \rightarrow \T$.
In Section~\ref{SS: change},
we describe the change of left key 
when we apply the reverse insertion
to a tableau.
In Section~\ref{SS: Proof of Theorem},
we prove Theorem~\ref{thm:insertionbijection}.

\subsection{Describing the reverse row insertion}
\label{SS: reverse insertion}

Huang, Shimozono and Yu~\cite{HSY}
introduced a row analogue of the Hecke column insertion~\cite{BKSTY}.
Notice that in~\cite{HSY}, the authors described the algorithm on decreasing tableaux.
For the purpose of this paper, we reverse all the comparisons
and describe the algorithm on increasing tableaux.
We start with the following definitions. 
For a tableau $P$,
we let $P_{r \downarrow}$
denote the tableau obtained by removing the first $r-1$ rows of $P$.

\begin{defn}[\cite{HSY}]
A value $x$ is \definition{ejectable} in $P$ if $x$ occurs in the first row of $P$ and either
$x+1$ is not in the first row of $P$, or $x+1$ is in the first row of $P$ and $x+1$ is ejectable in $P_{2 \downarrow}$.
\end{defn}

\begin{defn} Let $(r, c)$ be a cell in an increasing tableau $P$.
The \definition{bumping path} of $(r, c)$ in $P$
is the following sequence of cells in $P$: 
$(r, c_r), (r-1, c_{r-1}), \dots, (1, c_1)$
defined recursively. 
First, $c_r := c$.
Then $c_i$ is the largest such that $P(i, c_i) < P(i+1, c_{i+1})$.
\end{defn}

Clearly, a bumping path 
$(r, c_r), (r-1, c_{r-1}), \dots, (1, c_1)$
satisfies $c_r \leq \cdots \leq c_1$.

We say a cell $(r,c)$ in a tableau $P$ is an
\definition{outer cell} in $P$ if 
neither $(r+1,c)$ nor $(r, c+1)$ is in $P$.
The reverse row insertion algorithm takes $P$, $(r,c)$, and $\alpha$ as input
where $P$ is an increasing tableau, $(r,c)$ is an outer cell in $P$ and $\alpha = 0$ or $1$.
The output will be an increasing tableau $P'$ and 
a positive integer $m$. 
The input $\alpha$  indicates whether $P$ ``loses the cell $(r,c)$'':

\begin{align*}
    \shape(P') = \begin{cases} \shape(P) & \text{if $\alpha=0$} \\
    \shape(P) - \{(r,c)\} & \text{if $\alpha=1$.}
    \end{cases}
\end{align*}

Let $(r, c_r), (r-1, c_{r-1}), \dots, (1, c_1)$
be the bumping path starting at $(r,c)$ in $P$
and let $m_i = P(i, c_i)$.
The output value $m$ is by definition the value $m_1$.
The output tableau $P'$ will only differ from $P$ along the bumping path. It is enough to describe whether each $m_i$ on the bumping path gets replaced, and if so, by what value. This decision is determined iteratively by decreasing $i$ based on the values $m_i$ and $m_{i+1}$, the $i$-th row of $P$, the subtableau $P'_{>i}$, and a status indicator $\alpha_{i+1}\in \{0,1\}$. The $i$-th iteration updates the $i$-th row of $P$ (which becomes the $i$-th row of $P'$) and produces $\alpha_i\in \{0,1\}$.

Let $P'$ be a working tableau which is initialized to $P$.
In the initialization step, if $\alpha = 1$, remove from $P'$ the removable cell in row $r$ and its contents $m_r$ and set $\alpha_r=1$ and $i=r-1$.
If $\alpha=0$ set $m_{r+1} = \infty$, $\alpha_{r+1}=0$ and $i=r$.

The algorithm enters a loop. If $i = 0$ the algorithm terminates and the current 
tableau $P'$ is the output tableau. Now assume $i \geq 1$.
Let $R$ be the set consisting of numbers in row $i$ of the current tableau $P'$ (or equivalently $P$, since $P$ and $P'$ differ only in rows of index greater than $i$). By definition $m_i\in R$.
There are several cases for one iteration and each case has a nickname. 
\begin{enumerate}
\item[$\bullet$] \textbf{Dummy (Case D)}: If $m_i + 1 \in R$ (which implies $m_{i+1}=m_i + 1$) do not change the $i$-th row and set $\alpha_i = \alpha_{i+1}$.

\item[$\bullet$] \textbf{Direct Replacement (Case DR)}: Otherwise if $\alpha_{i+1}=1$ and $m_{i+1}\notin R$, replace $m_i$ by $m_{i+1}$ in row $i$ of $P'$ and set $\alpha_i=1$.
\end{enumerate}

Suppose neither of the two above cases hold. Find the largest ejectable entry $x$ in $P'_{i+1 \downarrow}$ such that $m_{i} < x < m_{i+1}$.
\begin{enumerate}
\item[$\bullet$] \textbf{Indirect Replacement (Case IR)}:
Suppose $x$ exists.
Replace $m_i$ by $x$ in row $i$ of $P'$ and set $\alpha_i=1$.
\item[$\bullet$] \textbf{No Replacement (Case NR)}: 
Suppose $x$ does not exist. Do not change the $i$-th row and set $\alpha_i=0$.
\end{enumerate}
Now decrement $i$ and go to the top of the loop.

\begin{exa}
Let $P$ be the increasing tableau below.
\[
P = 
\raisebox{0.9cm}{
\begin{ytableau}
1 & 2 & *(yellow)3 & 5\\
2 & *(yellow)5 & 6\\
3 & *(yellow)6\\
6 & *(yellow)7\\
8
\end{ytableau}}
\]
Let $(r,c)$ be $(4,2)$ and $\alpha = 0$.
First, the algorithm finds the bumping path starting at $(4,2)$,
which is highlighted by yellow. 
Thus, the output number $m$ is 3.
To compute the output tableau, 
it first initializes: $i = 4, \alpha_5 = 0$ and $m_5 = \infty$.
Now We trace the iterations of this algorithm.
\begin{itemize}
\item Iteration when $i = 4$: We have $m_4 = 7$.
Since $m_4 + 1 \notin R$,
we are not in Case D. 
Since $\alpha_5 = 0$,
we are not in Case DR.
Then find the largest $x$ 
that is ejectable in $P_{5 \downarrow}'$
and $7 < x < \infty$.
We have $x = 8$,
so we are in \textbf{Case IR}.
Replace the 7 in row 4 by 8
and set $\alpha_4 = 1, i = 3$.

\item Iteration when $i = 3$: We have $m_3 = 6$.
Since $m_3 + 1 \notin R$,
we are not in Case D. 
Since $\alpha_4 = 1$ and $m_4 \notin R$,
we are in \textbf{Case DR}.
Replace the 6 in row 3 by 7
and set $\alpha_3 = 1, i = 2$.

\item Iteration when $i = 2$: We have $m_2 = 5$.
Since $m_2 + 1 = 6$ is in $R$,
we are in \textbf{Case D}. 
We do not change row 2 
and set $\alpha_3 = 1, i = 1$.

\item Iteration when $i = 1$: We have $m_1 = 3$.
Since $m_1 + 1 \notin R$,
we are not in Case D. 
Since $m_2 \in R$,
we are not in Case DR.
Then find the largest $x$ 
that is ejectable in $P_{2 \downarrow}'$
and $3 < x < 5$.
Such $x$ does not exist,
so we are in \textbf{Case NR}.
We do not change row 1.
\end{itemize}
The tableau $P'$ 
after each iteration is depicted as follows:
\[
\begin{ytableau}
1 & 2 & 3 & 5\\
2 & 5 & 6\\
3 & 6\\
6 & *(green)8\\
8
\end{ytableau}
\quad\quad
\begin{ytableau}
1 & 2 & 3 & 5\\
2 & 5 & 6\\
3 & *(green)7\\
6 & 8\\
8
\end{ytableau}
\quad\quad
\begin{ytableau}
1 & 2 & 3 & 5\\
2 & 5 & 6\\
3 & 7\\
6 & 8\\
8
\end{ytableau}
\quad\quad
\begin{ytableau}
1 & 2 & 3 & 5\\
2 & 5 & 6\\
3 & 7\\
6 & 8\\
8
\end{ytableau}
\]
\[
\small\textrm{After 1 iteration}\quad\:
\small\textrm{After 2 iterations}\quad\:
\small\textrm{After 3 iterations}\quad\:
\small\textrm{After 4 iterations}
\]
\end{exa}

Huang, Shimozono and 
Yu also defined the insertion algorithm, which
takes an increasing tableau $P'$ 
and a positive integer $m$ as input,
producing a triple $(P, (r,c), \alpha)$.
By Theorem~5.3 of~\cite{HSY},
the reverse insertion algorithm and insertion are inverses of each other. 
We deduce the following property of 
the insertion algorithm.
\begin{lemma}
\label{L: restriction to smaller number (one insertion)}
Let $P$ be an increasing tableau
without the number $N$.
Perform row insertion on $(P, m)$
and get $(P', (r,c), \alpha)$.
Let $p$ (resp. $p'$) be the increasing 
tableau we get after ignoring all numbers larger than $N$ in $P$ 
(resp. $p'$).
If $m \geq N$,
then $p = p'$.
If $m < N$,
then $p'$ is the tableau we get 
after performing row insertion on 
$(p, m)$
\end{lemma}
\begin{proof}
We know $(P, m)$
is obtained by reverse row insertion 
on $(P', (r,c),\alpha)$.
If $m \geq N$,
then the whole bumping path does not involve numbers smaller than $N$.
Since $P$ is obtained from $P'$
by changing numbers in the bumping path,
we know $p = p'$.

Now suppose $m < N$
and consider the reverse insertion on
$(P', (r,c),\alpha)$.
Let $i$ be the largest such that
$m_i < N$.
Assume this value is at $(i, j)$ of $P$.
The algorithm behaves the same in the first $i$ rows of $P$ as computing 
the reverse insertion on $(p', (i,j), \alpha_{i+1})$.
Thus, reverse insertion on $(p', (i,j), \alpha_{i+1})$
yields $(p,m)$,
so $p'$
is the tableau we get after performing insertion on $(p,m)$.

\end{proof}

Now we may describe the bijection $\Psi$ 
from $\T$ to $\C$ recursively. 
Given $(P,Q) \in \T$.
If $P$ and $Q$ are both the empty tableaux,
then $\Psi(P,Q)$ is the pair of two empty words.
Otherwise, find the smallest number in $Q$
and break ties by picking the rightmost. 
Say it is the number $q$ in cell $(r,c)$.
Consider two cases.
\begin{itemize}
    \item If $q$ is the only number in $(r,c)$
    of $Q$, then we remove this cell
    and let $Q'$ be the resulting tableau.
    Perform the reverse insertion on $P$
    with input $(r,c)$ and $\alpha = 1$.
    Let $P', m$ be the output. 
    \item If $q$ is not the only number in $(r,c)$,
    we remove $q$ from this cell and 
    let $Q'$ be the resulting tableau.
    Perform the reverse insertion on $P$
    with input $(r,c)$ and $\alpha = 0$.
    Let $P', m$ be the output. 
\end{itemize}

For $(P', Q') \in \T$,
we may apply the map $\Psi$ recursively and get
$\Psi(P', Q') = (a', i')$.
Then $\Psi(P,Q) := (m \circ a', q \circ i')$
where $\circ$ represents concatenation of words. 

\begin{thm}\cite[Corollary 5.9]{HSY}
The map $\Psi$ is a bijection from
$\T_w$ to $\C_w$
that preserves the weight of the second entry.
In other words, $\Psi$ is bijective and 
if $\Psi(P,Q) = (a,i)$,
then $[\rev(\word(P))]_H= [a]_H$
and $\wt(Q) = \wt(i)$.
\end{thm}

\begin{exa}
Consider the following $(P,Q) \in \T$:
$$
P = 
\begin{ytableau}
1 & 2 \cr
3
\end{ytableau}\quad\quad
Q = 
\begin{ytableau}
3 & 21 \cr
21
\end{ytableau}\:.
$$
We would like to compute $\Psi(P,Q)$.
First, we find the smallest number in $Q$
and break ties by picking the rightmost.
That is the $1$ in cell $(1,2)$.
Thus, we perform reverse insertion on $P$
with input $(r,c) = (1,2)$ and $\alpha = 0$.
We get the output number $2$ and the resulting $(P',Q')$ is:
$$
P' = 
\begin{ytableau}
1 & 3 \cr
3
\end{ytableau}\quad\quad
Q' = 
\begin{ytableau}
3 & 2 \cr
21
\end{ytableau}\:.
$$

Recursively, we can compute $\Psi(P', Q') = (1313,1223)$.
Thus, $\Psi(P,Q)$ is $(21313, 11223)$.
The word $\rev(\word(P)) = 213$ and $21313$ are Hecke words 
of the same permutation. 
The RSVT $Q$ and the word $11223$ have the same weight. 
\end{exa}

Finally, we can prove Lemma~\ref{L: restriction to smaller number}
stated in \S\ref{S: Background}.
\begin{proof}[Proof of Lemma~\ref{L: restriction to smaller number}]
By induction on the length of $A$.
The inductive step is finished by Lemma~\ref{L: restriction to smaller number (one insertion)}.
\end{proof}

\subsection{Change of left key under reverse insertion}
\label{SS: change}

Shimozono and Yu~\cite{SY} studied how the left key
of an RSVT changes if we remove the rightmost appearance
of the smallest number.
\begin{lemma}[{\cite[Lemma 4.18 and Lemma 4.19]{SY}}]
\label{thm:reverseinsertkey Q}
Let $Q$ be an RSVT. 
Let $\min(Q) = i$ and let $(r,c)$
be the rightmost cell with $i \in Q(r,c)$.
Let $Q'$ be the RSVT obtained by removing $i$
from $Q$ (also remove the cell $(r,c)$ if $i$
is the only number in it).
If $i$ is not the only number in $Q$,
then clearly $K_-(L(Q)) = K_-(L(Q'))$.
Otherwise, $K_-(L(Q))$ and $K_-(L(Q'))$
only differ at column $c$:
$K_-(L(Q))_c = K_-(L(Q'))_c \sqcup y$
where
$$
y = \begin{cases}
i & \text{if $c = 1$,} \\
\min(K_-(L(Q'))_{c-1} - K_-(L(Q'))_{c}) &\text{if $c > 1$.}
\end{cases}
$$
\end{lemma}

We are going to derive the analogue of this lemma
for increasing tableau.
Take an increasing tableau $P$.
Fix an outer cell $(R,C)$ of $P$ and $\alpha = 0$ or $1$ throughout this section. 
Let $(P', m)$ be the result of reverse insertion
with input $(P, (R,C), \alpha)$.
We would like to show the following:
\begin{thm}
\label{thm:reverseinsertkey}
If $\alpha = 0$,
$K_-(P) = K_-(P')$.
Otherwise, 
$K_-(P)$ and $K_-(P')$ only differ
at column $C$:
$$K_-(P')_C = K_-(P)_C - \min(K_-(P)_C \setminus K_-(P)_{C+1}).$$
\end{thm}

The rest of this section
aims to prove Theorem~\ref{thm:reverseinsertkey}.
We define the \definition{block} in column $c$
as the set of $(r,c)$ such that $P(r,c) \neq P'(r,c)$.
Recall the definition of bumping path
from section \S\ref{S: Background}.
Let $(R, c_R), (R-1, c_{R-1}), \cdots, (1, c_1)$
be the bumping path of $P$ that starts at $(R,C)$.
Since the algorithm only changes entries
along the bumping path, we know
each block consists of cells in the bumping path.
Also recall that $c_R \leq c_{R-1} \leq \cdots \leq c_1$,
so the block in column $c$, if non-empty, must be 
$$\{(r_1, c), (r_1 + 1, c), \cdots, (r_2 - 1, c), (r_2, c)\}$$
for some $r_1 \leq r_2$.
Moreover, we know 
the reverse insertion performs direct replacement 
in row $r_1, \dots, r_2 - 1$.
Thus, for any $ r_1 \leq r < r_2$, 
$P'(r,c) = P(r + 1,c)$.
If $\alpha = 1$ and $c = C$,
$P'(r_2, c)$ is undefined. 
This block is also called the \definition{deletion block.}
Otherwise, the reverse insertion
performs direct replacement 
or indirect replacement on row $r_2$.
Such a block is also called an \definition{insertion block.}

\begin{exa}
Below is an example of the reverse-insertion algorithm applied on the cell $(9,1)$ of $P$ with $\alpha = 1$. We color the cells along the bumping path. The color is determined by the case of the algorithm in that row: red (DR), orange (IR), green (D), blue (NR), and pink (the initial removal).
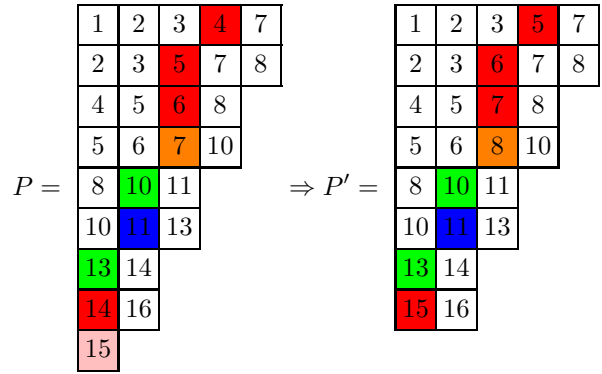
\[
P = 
\raisebox{2cm}{
\begin{ytableau}
1 & 2 & 3 & *(red)4 & 7\\
2 & 3 & *(red)5 & 7 & 8\\
4 & 5 & *(red)6 & 8\\
5 & 6 & *(orange)7 & 10\\
8 & *(green)10 & 11\\
10 & *(blue)11 & 13\\
*(green)13 & 14\\
*(red)14 & 16\\
*(pink)15
\end{ytableau}}
\Rightarrow
P' =
\raisebox{2cm}{
\begin{ytableau}
1 & 2 & 3 & *(red)5 & 7\\
2 & 3 & *(red)6 & 7 & 8\\
4 & 5 & *(red)7 & 8\\
5 & 6 & *(orange)8 & 10\\
8 & *(green)10 & 11\\
10 & *(blue)11 & 13\\
*(green)13 & 14\\
*(red)15 & 16\\
\end{ytableau}}
\]

Here, the deletion block is $\begin{ytableau} *(red) 14\\ *(pink) 15\end{ytableau}$ in the first column, and the insertion blocks are $\raisebox{0.4cm}{\begin{ytableau} *(red) 5\\ *(red) 6\\ *(orange) 7\end{ytableau}}$ in the third column, and $\raisebox{-0.1cm}{\begin{ytableau} *(red) 4\end{ytableau}}$ in the fourth column. Notice that outside of these blocks, $P$ and $P'$ are the same, and each block only changes by adding a number and removing a number from the block.
\end{exa}

We now view the effect of reverse insertion 
as changing the block in each column 
of $P$ to obtain $P'$.
We set $P^{(C-1)} = P$.
For $c \geq C$, let $P^{(c)}$
be the tableau obtained from $P^{c-1}$
by changing the block in column $c$ if it exists,
or $P^{(c)} = P^{(c-1)}$ if 
there is no block in column $c$.
Then $P' = P^{(m)}$ for some $m$ large enough.
To prove Theorem~\ref{thm:reverseinsertkey},
we just need to understand how $P^{(c)}$
differs from $P^{(c-1)}$.
First, we study the effect
of changing an insertion block.

\begin{lemma}
\label{lemma:insertionblock1}
Suppose the block in column $c$
is an insertion block and 
let $x$ be the value inserted 
to this block.  
Then $c>1$ and 
$|(P_c)_{<x}| = |(P'_{c-1})_{<x}|$.
\end{lemma}
\begin{proof}
Let $(r,c)$ be the bottom-most cell of this insertion block, so $P(r,c) = x$ 
and $P'(r,c) > x$.
Moreover, $P'(r, c-1) = P(r,c-1) < x$.
We know $x$ was inserted into $(r,c)$ through a direct or indirect replacement, so $P(r+1, c') = x$ for some $c'$. 
We consider two cases:
\begin{itemize}
\item Suppose the cell $(r+1, c)$ 
does not exist in $P$, 
then trivially $c' < c$
and $|(P_c)_{<x}| =r$.
If $(r+1, c-1)$ does not exist in $P$ either,
then $|(P_{c-1})_{<x}| =r$.
Otherwise, 
$P'(r+1, c-1) \geq P(r+1, c-1) \geq x$,
so $|(P'_{c-1})_{<x}| =r$.
\item Suppose the cell $(r+1, c)$ 
exists in $P$.
By $P'$ is increasing,
$x < P'(r+1, c) = P(r+1, c)$,
so $|(P_c)_{<x}| =r$.
Then $c' < c$ and
$x \leq P'(r+1, c-1)$,
so $|(P'_{c-1})_{<x}| =r$. \qedhere
\end{itemize}
\end{proof}

\begin{lemma}
\label{lemma:insertionblock3}
Suppose the block in column $c$
is an insertion block and 
let $x$ be the value inserted 
to this block. Then the insertion block removes the number $\min\{t\in P_c: x\sqcup P_c-t \preceq P_{c+1}\}$.
\end{lemma}
\begin{proof}
Let $z$ be the number removed by the insertion block
in column $c$.
We have $P'_c 
= x \sqcup P_c - z$
and $P'_c \preceq P_{c+1}$.
Assume $P(r,c) = z$. 
By definition of the insertion block, the algorithm cannot perform another direct replacement step in column $c$ in row $r-1$. Thus, either $r=1$ or $z \geq P_{c+1}(r-1)$.  In either case, for any $ i \in [r-1]$, we cannot have $x\sqcup P_c - P_c(i) \preceq P_{c+1}$.
Our claim is therefore proved. 
\end{proof}

We now understand the effect of changing an insertion block. 

\begin{lemma}
\label{lemma:insertionblock2}
Suppose the block in column $c$
is an insertion block,
then $K_-(P^{(c)}) = K_-(P^{(c-1)})$.
\end{lemma}
\begin{proof}
Say the insertion block gains $x$ and loses $y$,
so $P^{(c)}_{c} = x \sqcup P^{(c-1)}_c - y$.
By Corollary~\ref{C: leftkeycalc}, $K_-(P^{(c-1)})_i = P^{(c-1)}_{1,i}$ and $K_-(P^{(c)})_i = P^{(c)}_{1,i}$. Since $P^{(c-1)}$ and $P^{(c)}$ only differ in column $c$,
we know $K_-(P^{(c-1)})_i = K_-(P^{(c)})_i$
if $i < c$.

By Lemma~\ref{lemma:insertionblock1},
$c > 1$.
Notice that $P^{(c-1)}_{c-1} 
= P^{(c)}_{c-1} = P'_{c-1}$,
$P^{(c-1)}_{c+1} 
= P^{(c)}_{c+1} = P_{c+1}$,
$P^{(c-1)}_c = P_c$, 
and $P^{(c)}_c = P'_c$.
We set $U = P'_{c-1}$, 
$T = P_c$ and $S = P_{c+1}$.
By Lemma~\ref{lemma:insertionblock1}
and Lemma~\ref{lemma:insertionblock3},
$U, T, S$ and $x, y$ satisfy 
the conditions of 
Lemma~\ref{lemma:insert7}.
Thus,
$P'_{c-1} \lhd P_c = P'_{c-1} \lhd P'_c$
and 
$$
P'_{c-1} \lhd P_c \lhd S' 
= P'_{c-1} \lhd P'_c \lhd S'.
$$
for any $S' \subseteq P_{c+1}$.
The first equation implies 
$P'_{1, c} = P_{1, c}$.
The second equation implies
$P'_{1, i} = P_{1, i}$
if $i > c$.
\end{proof}

Next, we study
the deletion block.
\begin{lemma}
\label{lemma:deletionblock2}
Suppose the block
in column $c$ is a deletion block. 
Then this block loses $\min\{t\in P_c: P_c-t\preceq P_{c+1}\}$.
\end{lemma}
\begin{proof}
By a similar proof of Lemma~\ref{lemma:insertionblock3}.
\end{proof}

\begin{lemma}
\label{lemma:deletionblock3}
If a deletion block occurs in column $c$, then $K_-(P)$
and $K_-(P^{(c)})$ agree except in column $c$,
where
\begin{align}
\label{EQ: deletionblock3}
K_-(P^{(c)})_c = K_-(P)_c - \min(K_-(P)_c \setminus 
K_-(P)_{c+1})    
\end{align}

\end{lemma}
\begin{proof}

Recall that $P^{(c)}$ is obtained from $P$
by changing the deletion block. 
Assume the deletion block loses $y$, so $P^{(c)}_{c} = P_{c} -y$.  
By Lemma~\ref{lemma:deletionblock2}, $y= \min\{t\in P_c: P_c-t \preceq P_{c+1}\}$.

By Corollary~\ref{C: leftkeycalc}, 
column $i$ of $K_-(P)$ and 
$K_-(P^{(c)})$ are $P_{1,i}$ 
and $P^{(c)}_{1,i}$ respectively. 
They clearly agree if $i < c$.
Suppose $i > c$.
To show $P_{1,i} = P^{(c)}_{1, i}$,
it is enough to check $P_{c, i} =  P^{(c)}_{1, i}$, 
which is
$P_c \lhd P_{c+1,i} = (P_c - y) \lhd P_{c+1, i}$.
Since $P_{c+1, i} \subseteq P_{c+1}$,
this equation is implied by 
the third statement in Lemma~\ref{lemma:deletion7}.

Finally, we show 
$ P^{(c)}_{i,c} = P_{i,c} - \min(P_{i,c } \setminus P_{i,c+1})$ 
for all $i \in[c]$ by induction on $i$ from $i=c$ to $i=1$.  
Setting $i = 1$ would imply~(\ref{EQ: deletionblock3}). 
For the base case $i=c$, 
we know $P_{c,c}^{(c)} = P_{c,c} - y$.
It remains to check $y$ is the smallest
number in $P_c$ that is not in $P_{c, c+1}$.
This is done by the first two statements in \ref{lemma:deletion7}.
For the inductive step, assume $i < c$.  We have 
$$P^{(c)}_{i,c} = P_i \triangleleft P^{(c)}_{i+1,c} = P_i \triangleleft(P_{i+1,c}   - \min(P_{i+1,c} \setminus P_{i+1,c+1})) = P_{i,c}- \min(P_{i,c} \setminus P_{i,c+1}),$$ where the last
two steps are given by the inductive hypothesis and Lemma \ref{lemma:subsetkeycalc} respectively.
\end{proof}

Now we have understood how $K_-(P^{(c)})$
differs from $K_-(P^{(c-1)})$.

\begin{proof}[Proof of Theorem~\ref{thm:reverseinsertkey}]
Consider the sequence of increasing tableaux:
$$
P = P^{(C-1)}, P^{(C)}, P^{(C+1)}, \cdots, P^{(m-1)}, P^{(m)} = P',
$$
If $\alpha = 0$,
then there is no deletion block. 
By Lemma~\ref{lemma:insertionblock2},
$K_-(P) = K_-(P^{(C)}) = \cdots = K_-(P')$.
If $\alpha = 1$,
then the block in column $c$ is a
deletion block and all other blocks are insertion blocks.
By Lemma~\ref{lemma:insertionblock2},
$K_-(P^{(C)}) = K_-(P^{(C+1)}) = \cdots = K_-(P')$.
Then by Lemma~\ref{lemma:deletionblock3},
$K_-(P')$ agrees with $K_-(P)$ except in column $c$
where $K_-(P')_c = K_-(P)_c - \min(K_-(P)_c \setminus 
K_-(P)_{c+1})$.
\end{proof}

\subsection{Proof of Theorem~\ref{thm:insertionbijection}}
\label{SS: Proof of Theorem}
In this section, we prove Theorem~\ref{thm:insertionbijection}
using similar idea as~\cite[Theorem~4.2]{SY}.

\begin{proof}[Proof of Theorem~\ref{thm:insertionbijection}]
    Since $\Psi$ is a bijection from $\T$ to $\C$, it suffices to show that for any $(P,Q)\in \T^b$, $\Psi(P,Q)\in \C^b$ and for any $(a,i)\in \C^b$, $\Psi^{-1}(a,i) \in \T^b$.

    Take $(P,Q)\in \T^b$ where $Q$ has at least one cell. Let $i_1 = \min(Q)$. Let $(r,c)$ be the rightmost cell in $Q$ with $i_1 \in Q(r,c)$. Let $\alpha = 1$ if $i_1$ is the only number in $Q(r,c)$ and $\alpha = 0$ otherwise. 
    Let $(P', a_1)$ be the output of reverse insertion on $(P, (r,c), \alpha)$.
    Let $Q'$ be the tableau obtained from $Q$
    by removing $i_1$ from $Q(r,c)$,
    and also remove the cell $(r,c)$ if $\alpha = 1$.
    Due to the recursive definition of $\Psi$, 
    $\Psi(P,Q) = (a_1 \circ a', i_1 \circ i')$
    where $(a', i') = \Psi(P', Q')$.
    Notice that $(a,i) \in \C^b$
    is equivalent to $a_1 \leq i_1$
    and $(a', i') \in \C^b$.
    Inductively, we may assume $(a', i') \in \C^b$ if and only if $(P', Q') \in \T^b$.
    Thus, it suffices to show the following two statements:
    \begin{itemize}
    \item If $(P,Q) \in \T^b$,
    then $a_1 \geq i_1$ and $(P', Q') \in \T^b$.
    \item If $a_1 \geq i_1$ and $(P', Q') \in \T^b$,
    then $(P,Q) \in \T^b$.
    \end{itemize}

    We start with the first statement. 
    To see $a_1\geq i_1$, observe that  $i_1 = \min(Q)$ implies $i_1\leq \min(L(Q))$.  By the definition of $K_-(\cdot)$ on RSSYT, $\min(L(Q)) \leq \min(K_-(L(Q)))$, which in turn is at most $\min(K_-(P))  \leq \min(P) \leq a_1$.   
    Next we show $(P',Q')\in \T^b$.  If $\alpha = 0$, then Theorem \ref{thm:reverseinsertkey Q} and Theorem \ref{thm:reverseinsertkey} yield 
    $K_-(L(Q')) = K_-(L(Q))  \leq K_-(P) = K_-(P')$.
    Now suppose $\alpha = 1$.
    By Theorem \ref{thm:reverseinsertkey} the only column where $K_-(P)$ and $K_-(P')$ differ is column $c$ and $(K_-(P'))_c = (K_-(P))_c - x$, where $x =  \min((K_-(P))_c\setminus (K_-(P))_{c+1} )$. 
    By Theorem \ref{thm:reverseinsertkey Q},
    the only column where $K_-(L(Q))$ and $K_-(L(Q'))$ differ is column $c$ and $(K_-(Q'))_c = (K_-(Q))_c - y$ for some $y$.
    It remains to check for each $j \in [|K_-(P')_c|]$, $K_-(P')_c(j) \geq K_-(Q')_c(j)$.

    For $j$ such that $(K_-(P))_c(j)\geq x$, $$(K_-(P'))_c(j) = (K_-(P))_c(j+1) \geq (K_-(L(Q)))_c(j+1) \geq (K_-(L(Q')))_c(j).$$ 

For $j$ such that $(K_-(P))(j)< x$, $$(K_-(P'))_c(j) =  (K_-(P'))_{c+1}(j) \geq (K_-(L(Q')))_{c+1}(j) \geq   (K_-(L(Q')))_{c}(j).$$  Therefore, for all $j$, $(K_-(P'))_c(j) \geq  (K_-(L(Q')))_{c}(j)$, which shows that $K_-(P') \geq K_-(L(Q'))$.

We now show the second statement.  
If $\alpha = 0$, then Theorem \ref{thm:reverseinsertkey Q} and Theorem \ref{thm:reverseinsertkey} yield 
$K_-(L(Q')) = K_-(L(Q))  \leq K_-(P) = K_-(P')$.
    Otherwise, $\alpha = 1$. By Theorem \ref{thm:reverseinsertkey}, the only column where $K_-(P)$ and $K_-(P')$ differ is column $c$ and $(K_-(P))_c= (K_-(P'))_c\sqcup x$ for some $x$.  By Theorem \ref{thm:reverseinsertkey Q}, the only column where $K_-(L(Q))$ and $K_-(L(Q'))$ differ is column $c$, and $K_-(L(Q))_c= K_-(L(Q'))_c\sqcup y$,
    where $$
    y = \begin{cases}
    i_1 & \text{if $c = 1$,} \\
    \min(K_-(L(Q'))_{c-1} - K_-(L(Q'))_{c}) &\text{if $c > 1$.}
    \end{cases}
    $$
    
    It remains to check for each $j \in [|K_-(P)_c|]$, $K_-(P)_c(j) \geq K_-(Q)_c(j)$.
    For $j$ such that $(K_-(L(Q')))_c(j)> y$,
    $$(K_-(P))_c(j) \geq (K_-(P'))_c(j-1) \geq (K_-(L(Q')))_c(j-1) = (K_-(L(Q)))_c(j).$$
    Now assume $(K_-(L(Q')))_c(j) \leq y$.
    If $c = 1$, $y = i_1$ is the smallest
    number in $i$, so it is the smallest
    number in $Q$.
    We know $(K_-(L(Q')))_c(j) = y$ 
    and $y \leq \min(a) = \min(P) \leq (K_-(P))_c(j)$.
    If $c > 1$, we have
    $$
    K_-(P)_c(j) \geq K_-(P)_{c-1}(j) \geq K_-(L(Q))_{c-1}(j) = K_-(L(Q))_{c}(j)
    $$ 
    This completes the proof.
\end{proof}

 Finally, we mention a by-product of Theorem~\ref{thm:insertionbijection}: 
the Grothendieck-to-Lascoux expansion.
It was conjectured by~\cite{RY} 
and proved in~\cite{SY}.

\begin{cor}
\label{C: Gro-to-Las}
For a permutation $w \in S_n$,
\begin{equation*}
\fGb_w = \sum_{\substack{P \in \Inc \\ {\sf word}(P) = w^{-1}}} \beta^{|\wt(K_-(P))| - \ell(w)}\fLb_{\wt(K_-(P))}.
\end{equation*}
    \end{cor}
    \begin{proof}
        This follows from Theorem \ref{thm:insertionbijection}, Equation \ref{EQ: RSVT rule for Lascoux}, and Equation \ref{EQ: GrotDef}.
    \end{proof}

\section{Acknowledgements}
We thank Alex Yong for suggesting this problem.
We thank Brendon Rhoades and Alex Yong for
carefully reading an earlier version of this paper and 
giving many useful comments.
The first author was partially supported by NSF RTG DMS 1937241.

\appendix

\section{Example of Theorem \ref{thm:G times L}}
\label{appendix}
Suppose $n = 3$.
Say we would like to expand 
$\fLb_{\alpha} \times 
G^{(\beta)}_{w}(x_1, \cdots, x_n)$ 
into Lascoux polynomials
where $\alpha = (1,0,2)$ 
and $w$ has one-line notation
$12345876$.
We let 
$$
P_1 = 
\raisebox{0.1cm}{
\begin{ytableau}
1 & 4\\
3
\end{ytableau}},
$$
so $K_-(P_1) = (1,0,2)$.
We pick $N = 5$.
Now we are looking for increasing
tableau with at most $3$ rows such that
\begin{itemize}
\item Entries smaller than $5$
form the $P_1$ we picked, and
\item Entries of $\word(P)$
that are larger than $5$ form a Hecke word
of $w$.
\end{itemize}

The following are all increasing 
tableaux satisfying both conditions:
$$
\begin{ytableau}
1 & 4 & 6 & 7\\
3 \\
7
\end{ytableau}
\quad\quad
\begin{ytableau}
1 & 4 & 6 & 7\\
3 & 7\\
\end{ytableau}
\quad\quad
\begin{ytableau}
1 & 4 & 6\\
3 & 6 & 7\\
\end{ytableau}
\quad\quad
\begin{ytableau}
1 & 4 & 6\\
3 & 7\\
6
\end{ytableau}
\quad\quad
\begin{ytableau}
1 & 4 & 7\\
3 & 6\\
7
\end{ytableau}
\quad\quad
\begin{ytableau}
1 & 4\\
3 & 6\\
6 & 7
\end{ytableau}
$$

$$
\begin{ytableau}
1 & 4 & 6 & 7\\
3 & 7\\
7
\end{ytableau}
\quad\quad\quad\quad
\begin{ytableau}
1 & 4 & 6 & 7\\
3 & 6\\
7
\end{ytableau}
\quad\quad\quad\quad
\begin{ytableau}
1 & 4 & 6 & 7\\
3 & 7\\
6
\end{ytableau}
\quad\quad\quad\quad
\begin{ytableau}
1 & 4 & 6 & 7\\
3 & 6 & 7\\
\end{ytableau}
$$$$
\begin{ytableau}
1 & 4 & 6 \\
3 & 6\\
6 & 7
\end{ytableau}
\quad\quad\quad\quad
\begin{ytableau}
1 & 4 & 7 \\
3 & 6\\
6 & 7
\end{ytableau}
\quad\quad\quad\quad
\begin{ytableau}
1 & 4 & 6 \\
3 & 6 & 7\\
6
\end{ytableau}
\quad\quad\quad\quad
\begin{ytableau}
1 & 4 & 6 \\
3 & 6 & 7\\
7
\end{ytableau}
$$

$$
\begin{ytableau}
1 & 4 & 6 & 7\\
3 & 6 & 7\\
6
\end{ytableau}
\quad\quad\quad\quad
\begin{ytableau}
1 & 4 & 6 & 7\\
3 & 6 & 7\\
7
\end{ytableau}
\quad\quad\quad\quad
\begin{ytableau}
1 & 4 & 6 & 7\\
3 & 6 \\
6 & 7
\end{ytableau}
\quad\quad\quad\quad
\begin{ytableau}
1 & 4 & 6 \\
3 & 6 & 7\\
6 & 7
\end{ytableau}
$$

$$
\begin{ytableau}
1 & 4 & 6 & 7\\
3 & 6 & 7\\
6 & 7
\end{ytableau}
$$
After applying $K_-(\cdot)$, we obtain the keys:
$$
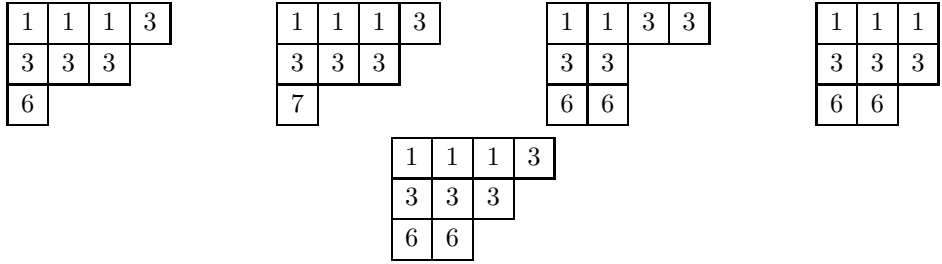
\begin{ytableau}
1 & 3 & 3 & 3\\
3 \\
7
\end{ytableau}
\quad\quad
\begin{ytableau}
1 & 1 & 3 & 3\\
3 & 3\\
\end{ytableau}
\quad\quad
\begin{ytableau}
1 & 1 & 1\\
3 & 3 & 3\\
\end{ytableau}
\quad\quad
\begin{ytableau}
1 & 3 & 3\\
3 & 6\\
6
\end{ytableau}
\quad\quad
\begin{ytableau}
1 & 1 & 3\\
3 & 3\\
7
\end{ytableau}
\quad\quad
\begin{ytableau}
1 & 1\\
3 & 3\\
6 & 6
\end{ytableau}
$$

$$
\begin{ytableau}
1 & 1 & 3 & 3\\
3 & 3\\
7
\end{ytableau}
\quad\quad\quad\quad
\begin{ytableau}
1 & 1 & 3 & 3\\
3 & 3\\
7
\end{ytableau}
\quad\quad\quad\quad
\begin{ytableau}
1 & 3 & 3 & 3\\
3 & 6\\
6
\end{ytableau}
\quad\quad\quad\quad
\begin{ytableau}
1 & 1 & 1 & 3\\
3 & 3 & 3\\
\end{ytableau}
$$$$
\begin{ytableau}
1 & 1 & 3 \\
3 & 3\\
6 & 6
\end{ytableau}
\quad\quad\quad\quad
\begin{ytableau}
1 & 1 & 3 \\
3 & 3\\
6 & 6
\end{ytableau}
\quad\quad\quad\quad
\begin{ytableau}
1 & 1 & 1\\
3 & 3 & 3\\
6
\end{ytableau}
\quad\quad\quad\quad
\begin{ytableau}
1 & 1 & 1 \\
3 & 3 & 3\\
7
\end{ytableau}
$$

$$
\begin{ytableau}
1 & 1 & 1 & 3\\
3 & 3 & 3\\
6
\end{ytableau}
\quad\quad\quad\quad
\begin{ytableau}
1 & 1 & 1 & 3\\
3 & 3 & 3\\
7
\end{ytableau}
\quad\quad\quad\quad
\begin{ytableau}
1 & 1 & 3 & 3\\
3 & 3 \\
6 & 6
\end{ytableau}
\quad\quad\quad\quad
\begin{ytableau}
1 & 1 & 1 \\
3 & 3 & 3\\
6 & 6
\end{ytableau}
$$

$$
\begin{ytableau}
1 & 1 & 1 & 3\\
3 & 3 & 3\\
6 & 6
\end{ytableau}
$$

We apply $\textrm{cap}_3$ to each key and get:
$$
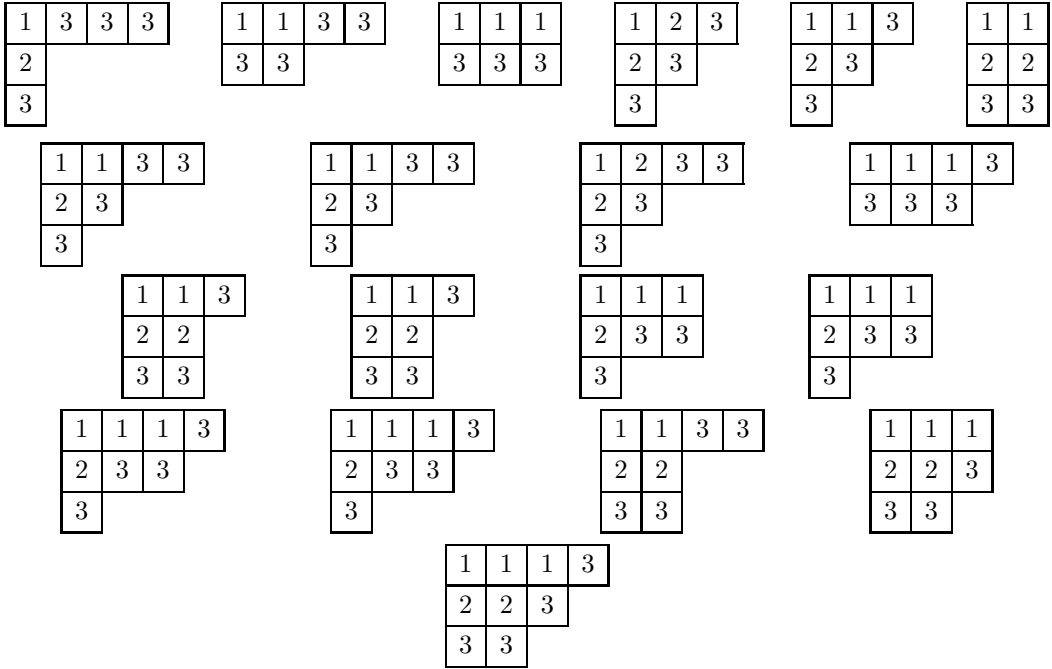
\begin{ytableau}
1 & 3 & 3 & 3\\
2 \\
3
\end{ytableau}
\quad\quad
\begin{ytableau}
1 & 1 & 3 & 3\\
3 & 3\\
\end{ytableau}
\quad\quad
\begin{ytableau}
1 & 1 & 1\\
3 & 3 & 3\\
\end{ytableau}
\quad\quad
\begin{ytableau}
1 & 2 & 3\\
2 & 3\\
3
\end{ytableau}
\quad\quad
\begin{ytableau}
1 & 1 & 3\\
2 & 3\\
3
\end{ytableau}
\quad\quad
\begin{ytableau}
1 & 1\\
2 & 2\\
3 & 3
\end{ytableau}
$$

$$
\begin{ytableau}
1 & 1 & 3 & 3\\
2 & 3\\
3
\end{ytableau}
\quad\quad\quad\quad
\begin{ytableau}
1 & 1 & 3 & 3\\
2 & 3\\
3
\end{ytableau}
\quad\quad\quad\quad
\begin{ytableau}
1 & 2 & 3 & 3\\
2 & 3\\
3
\end{ytableau}
\quad\quad\quad\quad
\begin{ytableau}
1 & 1 & 1 & 3\\
3 & 3 & 3\\
\end{ytableau}
$$$$
\begin{ytableau}
1 & 1 & 3 \\
2 & 2\\
3 & 3
\end{ytableau}
\quad\quad\quad\quad
\begin{ytableau}
1 & 1 & 3 \\
2 & 2\\
3 & 3
\end{ytableau}
\quad\quad\quad\quad
\begin{ytableau}
1 & 1 & 1\\
2 & 3 & 3\\
3
\end{ytableau}
\quad\quad\quad\quad
\begin{ytableau}
1 & 1 & 1 \\
2 & 3 & 3\\
3
\end{ytableau}
$$

$$
\begin{ytableau}
1 & 1 & 1 & 3\\
2 & 3 & 3\\
3
\end{ytableau}
\quad\quad\quad\quad
\begin{ytableau}
1 & 1 & 1 & 3\\
2 & 3 & 3\\
3
\end{ytableau}
\quad\quad\quad\quad
\begin{ytableau}
1 & 1 & 3 & 3\\
2 & 2 \\
3 & 3
\end{ytableau}
\quad\quad\quad\quad
\begin{ytableau}
1 & 1 & 1 \\
2 & 2 & 3\\
3 & 3
\end{ytableau}
$$

$$
\begin{ytableau}
1 & 1 & 1 & 3\\
2 & 2 & 3\\
3 & 3
\end{ytableau}
$$

Finally, we apply $\wt(\cdot)$ and get:
$$
(1,1,4),\quad\quad (2,0,4),\quad\quad(3,0,3),\quad\quad(1,2,3),\quad\quad(2,1,3),\quad\quad(2,2,2)
$$
$$
(2,1,4),\quad\quad\quad\quad (2,1,4), \quad\quad\quad\quad (1,2,4), \quad\quad\quad\quad (3,0,4)
$$
$$
(2,2,3), \quad\quad\quad\quad(2,2,3), \quad\quad\quad\quad(3,1,3), \quad\quad\quad\quad(3,1,3)
$$
$$
(3,1,4), \quad\quad\quad\quad(3,1,4), \quad\quad\quad\quad(2,2,4), \quad\quad\quad\quad(3,2,3)
$$
$$
(3,2,4)
$$

As a conclusion, we expand $\fLb_\alpha \times G^{(\beta)}_{w}(x_1, x_2, x_3)$ into:
\begin{align*}
& \fLb_{(1,1,4)} +  \fLb_{(2,0,4)} + \fLb_{(3,0,3)} +  \fLb_{(1,2,3)} +  \fLb_{(2,1,3)} + \fLb_{(2,2,2)}\\
+ &  2 \beta\fLb_{(2,1,4)} +  \beta\fLb_{(1,2,4)} + 2\beta\fLb_{(2,2,3)} +  2\beta\fLb_{(3,1,3)}\\
+ & 2\beta^2\fLb_{(3,1,4)} +  \beta^2\fLb_{(2,2,4)} + \beta^2\fLb_{(3,2,3)}
+ \beta^3 \fLb_{(3,2,4)}
\end{align*}

\bibliographystyle{alpha}
\bibliography{main.bbl}{}
\end{document}